\documentclass[11pt]{article}
\usepackage{amsfonts}
\usepackage{amsfonts,latexsym,amsmath,amscd,geometry}
\usepackage{amssymb}
\usepackage{latexsym}
\renewcommand{\theequation}{\thesection.\arabic{equation}}
\newcommand \nc{\newcommand}
\newtheorem{theorem}{Theorem}[section]
\newtheorem{lemma}[theorem]{Lemma}

\newtheorem{corollary}[theorem]{Corollary}
\newtheorem{definition}[theorem]{Definition}

\newtheorem{remark}[theorem]{Remark}
\renewcommand{\thetheorem}{\thesubsection.\arabic{theorem}}
\nc{\ba}{\begin{array}}\nc{\ea}{\end{array}}
\nc{\be}{\begin{eqnarray}}\nc{\ee}{\end{eqnarray}}
\nc{\beq}{\begin{equation}}\nc{\eeq}{\end{equation}}
\nc{\bex}{\begin{eqnarray*}}\nc{\eex}{\end{eqnarray*}}
\nc{\btm}{\begin{theorem}} \nc{\etm}{\end{theorem}}
\nc{\blm}{\begin{lemma}} \nc{\elm}{\end{lemma}}
\nc{\R}{\mathbb{R}} \nc{\va}{\varepsilon} \nc{\ls}{\limits}

\def\de{\Delta}
\def\pf{\noindent{\bf Proof.\quad}}\def\endpf{\hfill$\Box$}
\def\les{\lesssim}\def\u{\dot{u}}\def\di{\mbox{div\,}}

\begin{document}
\title{\bf Blow-up criterions of strong solutions to 3D compressible Navier-Stokes equations with vacuum}
\author{
\begin{tabular}{cc}
& H{\sc uanyao} W{\sc en},  \ \ \ C{\sc hangjiang} Z{\sc
hu}\thanks{Corresponding author.\ \ Email:
cjzhu@mail.ccnu.edu.cn} \\[4mm]
& The Hubei Key Laboratory of Mathematical Physics \\
& School of Mathematics and
Statistics \\
& Central China Normal University, Wuhan 430079, P.R. China\\
\end{tabular}
}
\date{}

\maketitle

\begin{abstract}
In the paper, we establish a blow-up criterion in terms of the integrability of the density for strong solutions to the
Cauchy problem of compressible isentropic Navier-Stokes equations in $\mathbb{R}^3$ with vacuum, under the assumptions
on the coefficients of viscosity: $\frac{29\mu}{3}>\lambda$. This extends the corresponding results in \cite{Huang-Li-Xin: Serrin, Sun-Wang-Zhang}
where a blow-up criterion in terms of the upper bound of the density was obtained under the condition $7\mu>\lambda$.
As a byproduct, the restriction $7\mu>\lambda$ in \cite{Fan-Jiang-Ou, Sun-Wang-Zhang 1} is relaxed to $\frac{29\mu}{3}>\lambda$
for the full compressible Navier-Stokes equations by giving a new proof of Lemma \ref{blow-up:le:5.1}. Besides, we get a blow-up criterion
in terms of the upper bound of the density and the temperature for strong solutions to the Cauchy problem of the full compressible
Navier-Stokes equations in $\mathbb{R}^3$. The appearance of vacuum could be allowed. This extends the corresponding results
in \cite{Sun-Wang-Zhang 1} where a blow-up criterion in terms of the upper bound of $(\rho,\frac{1}{\rho}, \theta)$
was obtained without vacuum. The effective viscous flux plays a very important role in the proofs.

\end{abstract}
\vspace{4mm}


\tableofcontents

\section {Introduction}
\setcounter{equation}{0} \setcounter{theorem}{0}
The full compressible Navier-Stokes
equations in $\mathbb{R}^N$ are written as follows:
 \be\label{full N-S}
\begin{cases}
\rho_t+\nabla \cdot (\rho u)=0, \\
(\rho  u)_t+\mathrm{div}(\rho u\otimes u)+\nabla
P=\mathrm{div}(\mathcal {T}),\\
(\rho E)_t+\mathrm{div}(\rho E u)+\mathrm{div}(P
u)=\mathrm{div}(\mathcal {T}
u)+\mathrm{div}(\kappa\nabla\theta).
\end{cases}
\ee Here $\mathcal {T}$ is the stress tensor, given by
$$\mathcal
{T}=\mu\left(\nabla  u+(\nabla u)^\prime\right)+\lambda
\di  u I_N,$$ where $I_N$ is a $N\times N$ unit matrix;
 $\rho=\rho(x,t)$, $u=u(x, t):\mathbb{R}^N\times(0,\infty)\rightarrow\mathbb{R}^N$, and $\theta=\theta(x,t)$
are unknown functions denoting  the density, velocity and absolute
temperature, respectively; $P$, $E$ and $\kappa$ denote
pressure, total energy and coefficient
of heat conduction, respectively, where $E=e+\frac{|u|^2}{2}$ ($e$ is the internal
energy), and $\kappa$ is a positive constant. Here, the state equations of $P$ and $e$ is of ideal polytropic gas type:
$$ P=a\rho\theta,\ e=C_0\theta,
$$
where $a$ and $C_0$ are two positive constants. $\mu$ and $\lambda$ are the coefficients of viscosity, which are assumed to be constants, satisfying the following physical restrictions:
$$
\mu>0,\ 2\mu+N\lambda\ge0.$$
For isentropic fluids, the compressible Navier-Stokes equations become
 \be\label{N-S-1}
\begin{cases}
\rho_t+\nabla \cdot (\rho u)=0, \\
(\rho  u)_t+\mathrm{div}(\rho u\otimes u)+\nabla
P=\mu\Delta u+(\mu+\lambda)\nabla\mathrm{div} u.
\end{cases}
\ee Here $P$ satisfies the equation of state of an ideal fluid:
$$
P=a\rho^\gamma,\ (a>0,\ \gamma>1).
$$
The compressible Navier-Stokes system is a well-known mathematical model which describes the
motion of compressible fluids (refer for instance to \cite{Lions} and references therein). There are so many known results on the well-posedness of solutions to (\ref{full N-S}) and (\ref{N-S-1}). In the absence of vacuum (vacuum means $\rho=0$), please refer for instance to \cite{Hoff1995, Hoff ARMA, Itaya, Jiang1, Kawohl,
Kazhikhov-Shelukhi, Matsumura-Nishida: Kyoto Un,
Matsumura-Nishida: CMP, Tani} and references therein.

We give a brief survey on the well-posedness of solutions to
(\ref{N-S-1}) and (\ref{full N-S}) with vacuum. First, for
(\ref{N-S-1}), there has been made great progress since Lions' work.
More precisely, the existence of global weak solutions to
(\ref{N-S-1}) with large initial data in $\mathbb{R}^N$ was first
obtained by Lions in \cite{Lions}, where $\gamma\ge\frac{3N}{N+2}$
for $N=2$ or 3. Feireisl et al in \cite{Feireisl2} extended Lions'
work to the case $\gamma>\frac{3}{2}$ for $N=3$. For solutions with
spherical symmetry, Jiang and Zhang in \cite{Jiang} relaxed the
restriction on $\gamma$ in \cite{Lions} to the case $\gamma>1$, and
got the global existence of the weak solutions for $N=2$ or $3$. On
the existence and regularity of weak solutions with density
connecting to vacuum continuously in 1D, please refer to \cite
{Luo}. During the pass two decades, Salvi, Choe, Kim and Jiang et al
made progress towards the local or global existence of strong
solutions with vacuum, see \cite{Cho-Choe-Kim, Choe-Kim:symmetric,
Fan-Jiang-Ni, salvi}. On the classical solutions, refer to
\cite{Cho-Kim: classical} for the local existence in three space
dimension, and refer to \cite{Huang-Li-Xin: classical} for global
existence with small initial energy in 3D, and refer to
\cite{Ding-Wen-Zhu} for global existence with large initial data in
1D. Secondly, for (\ref{full N-S}), the results on the global
existence of weak solutions can be referred to
\cite{Bresch-Desjardins, Feireisl-book}). More precisely, Feireisl
in \cite{Feireisl-book} got the global existence of {\em
variational} solutions in dimension $N\ge 2$. The temperature
equation in \cite{Feireisl-book} is satisfied only as an inequality
in the sense of distributions. Feireisl's work is the very first
attempt towards the global existence of weak solutions to the full
compressible Navier-Stokes equations in high dimensions. In order
that the equations are satisfied as equalities in the sense of
distribution, Bresch and Desjardins in \cite{Bresch-Desjardins}
proposed some different assumptions from \cite{Feireisl-book}, and
obtained the existence of global weak solutions to the full
compressible Navier-Stokes equations with large initial data and
density-dependent viscosities in $\mathbb{T}^3$ or $\mathbb{R}^3$.
On the regularities of the solutions to (\ref{full N-S}) when vacuum
is allowed, please refer to \cite{Cho-Kim; perfect gas} for the
local existence and uniqueness of strong solutions in bounded or
unbounded domains $\Omega\subseteq\mathbb{R}^3$, and refer to
\cite{Wen-Zhu 1} for the global existence and uniqueness of
classical solutions with large initial data in a bounded domain
$I\subseteq\mathbb{R}^1$, and refer to \cite{Wen-Zhu 2} for the the
global existence and uniqueness of spherically or cylindrically
symmetric classical solutions with large initial data in a bounded
domain $\Omega\subseteq\mathbb{R}^3$.

 It should be noted that one would not expect better regularities of the
solutions of (\ref{full N-S}) or (\ref{N-S-1}) in general because of Xin's results (\cite{Xin}) and Rozanova's results (\cite{Rona}).
It was proved that there is no global smooth
solution to the Cauchy problem of (\ref{full N-S}) or (\ref{N-S-1}), if the initial density is nontrivial compactly supported (\cite{Xin},
$N=1$ for (\ref{N-S-1}) and $N\ge1$ for (\ref{full N-S})) or the solutions are highly decreasing at infinity (\cite{Rona}, $N\ge3$ for (\ref{N-S-1})
and (\ref{full N-S})). In fact, a similar problem which is largely open for the incompressible Navier-Stokes equations in $\mathbb{R}^3$, i.e., whether the global smooth solutions exist or not, was proposed as one of the Millennium Prize Problems by Clay Mathematics Institute (CMI) (see \cite{Carlson}, 57-67: Charles L. Fefferman, Existence and Smoothness of the Naiver-Stokes Equation). These motivate us to find some possible blow-up criterions of regular solutions to (\ref{full N-S}) and (\ref{N-S-1}), especially
of strong solutions. Such a problem has been studied for the incompressible Euler equations by Beale-Kato-Majda in their poineering work \cite{BKM},
which showed that the $L^1_tL^\infty_x$-bound of vorticity $\nabla\times u$ alone controls the breakdown of smooth solutions. Later, Ponce \cite{Ponce}
rephrased the BKM-criterion in terms of the deformation tensor $T_{ij}=\partial_ju^i+\partial_iu^j$. Recently, some results on the blow-up criterions have been done for some related models, such as compressible liquid crystal system which is the one coupling compressible Navier-Stokes equations with heat flow of harmonic map, see for instance \cite{Huang-Wang-Wen, Huang-Wang-Wen:arma}.

Before stating our main result, We would like to give some notations which will be used throughout
the paper.
\section {Main results}
\setcounter{equation}{0} \setcounter{theorem}{0}
Before stating our main results, We would like to give some notations which will be used throughout
the paper.

\subsection {Notations}
\setcounter{equation}{0} \setcounter{theorem}{0}

\ \ \ \ \ (i)\ $\int_{\mathbb{R}^3} f =\int_{\mathbb{R}^3} f \,dx.$\\

(ii)\ For $1\le l\le \infty$, denote the $L^l$ spaces and the standard Sobolev spaces as follows:
$$L^l=L^l(\mathbb{R}^3),  \ D^{k,l}=\left\{ u\in L^1_{\rm{loc}}(\mathbb{R}^3): \|\nabla^k u \|_{L^l}<\infty\right\},$$
$$W^{k,l}=L^l\cap D^{k,l},  \ H^k=W^{k,2}, \ D^k=D^{k,2},$$
$$D_0^1=\Big\{u\in L^6: \ \|\nabla u\|_{L^2}<\infty\},$$
$$\|u\|_{D^{k,l}}=\|\nabla^k u\|_{L^l}.$$

(iii)\ For two $3\times 3$ matrices $E=(E_{ij}), F=(F_{ij})$, denote the scalar product between
$E$ and $F$ by
$$E:F=\sum_{i,j=1}^3 E_{ij} F_{ij}.$$

(iv)\ $G=(2\mu+\lambda)\mathrm{div}u-P$ is the effective viscous flux.\\

(v)\ $\dot{h}=h_t+u\cdot\nabla h$ denotes the material derivative.

\subsection {Compressible isentropic N-S: a blow-up criterion in terms of the integrability of the density}
\setcounter{equation}{0} \setcounter{theorem}{0}
\renewcommand{\theequation}{\thesubsection.\arabic{equation}}
The constant $a$ in the pressure function plays no roles in the analysis, we assume $a=1$ henceforth. If the solutions are regular enough (such as strong solutions), (\ref{N-S-1}) is equivalence to the following system which is very usefull in the proofs of the main theorems:
\be\label{N-S}
\begin{cases}
\rho_t+\nabla \cdot (\rho u)=0, \\
\rho  u_t+\rho u\cdot\nabla u+\nabla
P=\mu\Delta u+(\mu+\lambda)\nabla\mathrm{div} u,\ \mathrm{in}\ \mathbb{R}^3.
\end{cases}
\ee
System (\ref{N-S}) is supplemented with initial conditions
\be\label{initial}
(\rho, u)|_{t=0}=(\rho_0, u_0),\ x\in\mathbb{R}^3,
\ee
with
\be\label{boundary}
\rho(x,t)\rightarrow0,\ u(x,t)\rightarrow0\ \mathrm{as}\ |x|\rightarrow\infty,\ t\ge0.
\ee
We give the definition of strong solutions to (\ref{N-S}) throughout the paper.
\begin{definition}(Strong solutions) For $T>0$, $(\rho, u)$ is called a strong solution to the compressible Navier-Stokes equations (\ref{N-S})-(\ref{boundary}) in $\mathbb{R}^3\times [0,T]$, if for some $q\in (3, 6]$,
\bex &0\le \rho\in
C([0,T];W^{1,q}\cap H^1\cap L^1),\ \rho_t\in C([0,T];L^2\cap L^q),&\\
& u\in C([0,T];D^2\cap D^1_0)\cap L^2(0,T;D^{2,q}),\ u_t\in
L^2(0,T;D^1_0),\
 \sqrt{\rho} u_t\in
L^\infty(0,T;L^2),& \eex
and $(\rho,u)$ satisfies (\ref{N-S}) a.e. in $\mathbb{R}^3\times (0,T]$.
\end{definition}
Our main result for compressible isentropic Navier-Stokes equations is stated as follows:
\begin{theorem}\label{th:1.1}
Assume $\rho_0\geq0$,
 $\rho_0\in L^1\cap H^1\cap W^{1,q}$, for some $q\in(3,6]$, $u_0\in D^2\cap D_0^1$, and the
following compatibility conditions are satisfied:
 \beq\label{compatibility}
\mu\Delta u_0+(\mu+\lambda)\nabla\mathrm{div} u_0-\nabla P(\rho_0)=\sqrt{\rho}_0g,\ x\in\mathbb{R}^3,
\eeq for some $g\in L^2$. Let $(\rho, u)$ be a strong solution to (\ref{N-S})-(\ref{boundary}) in $\mathbb{R}^3\times[0, T]$. If $0<T^*<+\infty$ is the maximum time of existence of the strong solution,
then
\be\label{result}
\lim\sup\limits_{T\nearrow T^*}\|\rho\|_{L^\infty(0,T; L^{q_1})}=\infty,
\ee for some $1<q_1<\infty$ large enough, provided $\frac{29\mu}{3}>\lambda$.
\end{theorem}
\begin{remark}  Under the conditions of Theorem \ref{th:1.1},
the local existence of the strong solutions was obtained in \cite{Cho-Choe-Kim}. Thus, the assumption $T^*>0$ makes sense.
\end{remark}
In the presence of vacuum, before Theorem \ref{th:1.1}, there are several results on the blow-up criterions of strong solutions to (\ref{N-S}), refer for instance to \cite{Cho-Choe-Kim, Fan-Jiang, Huang-Li-Xin: Serrin, Huang-Li-Xin: blow up, Sun-Wang-Zhang}. More precisely, let $0<T^\star<+\infty$ is the maximum time of existence of strong solutions. Then the blow-up criterions can be summed as follows:

$\bullet$ Cho-Choe-Kim (\cite{Cho-Choe-Kim})
\be\label{cho criterion}
\lim\sup\limits_{t\nearrow T^\star}(\|\rho(t)\|_{H^1\cap
W^{1,q}}+\|u(t)\|_{D_0^1})=\infty,
\ee
for some $q\in(3,6]$;

$\bullet$ Fan-Jiang (\cite{Fan-Jiang})
\be\label{fan-jiang criterion}
\lim\sup\limits_{t\nearrow
T^\star}\left(\|\rho(t)\|_{L^\infty}+\int_0^t(\|\rho(s)\|_{W^{1,q}}+\|\nabla\rho(s)\|_{L^2}^4)\,ds\right)=\infty,\ee
for some $q\in(3,6]$, provided $7\mu>9\lambda$;

$\bullet$ Huang-Li-Xin (\cite{Huang-Li-Xin: blow up})
\be\label{Huang-Li-Xin criterion}
\lim\sup\limits_{t\nearrow T^\star}\int_0^t\|\frac{\nabla
u(s)+(\nabla u)^\prime(s)}{2}\|_{L^\infty}\,ds=\infty;
\ee

$\bullet$ Huang-Li-Xin (\cite{Huang-Li-Xin: Serrin}) (Serrin's criterion \cite{serrin})
\be\label{Huang-Xin criterion}
\lim\sup\limits_{t\nearrow
T^\star}\left(\|\rho\|_{L^\infty(0,t;L^\infty)}+\|\sqrt{\rho}u\|_{L^s(0,t;L^r)}\right)=\infty,
\ee
where $\frac{2}{s}+\frac{3}{r}\le1$, $3<r\le\infty$;

$\bullet$ Huang-Li-Xin (\cite{Huang-Li-Xin: Serrin}, for Cauchy problem), Sun-Wang-Zhang (\cite{Sun-Wang-Zhang}, for Cauchy problem and IBVP)
\be\label{Huang-Xin and Sun-Wang-Zhang criterion}
\lim\sup\limits_{t\nearrow
T^\star}\|\rho\|_{L^\infty(0,t;L^\infty)}=\infty,
\ee
provided $7\mu>\lambda$.\\

We introduce the main ideas of the proof of Theorem \ref{th:1.1}, some of which are inspired by
some of the arguments in \cite{Choe-Bum, Huang-Li-Xin: Serrin, Sun-Wang-Zhang, Wen-Yao-Zhu}.

(1) In \cite{Huang-Li-Xin: Serrin, Sun-Wang-Zhang}, to prove
(\ref{Huang-Xin and Sun-Wang-Zhang criterion}), the restriction
$7\mu>\lambda$ plays an important role in the analysis. In fact, the
condition $7\mu>\lambda$ is only used to get the upper bound of
$\int_{\mathbb{R}^3}\rho |u|^r$, for some $r>3$, so is it for
(\ref{Fan-Jiang-Ou criterion}) and (\ref{Sun-Wang-Zhang criterion}).
Here, we get the upper bound of $\int_{\mathbb{R}^3}\rho |u|^r$,
under the assumption $\frac{29\mu}{3}>\lambda$ (see Lemma
\ref{blow-up:le:5.1}), which as a byproduct of Lemma
\ref{blow-up:le:5.1} extends the results in \cite{Fan-Jiang-Ou,
Huang-Li-Xin: Serrin, Sun-Wang-Zhang, Sun-Wang-Zhang 1} (see Remark
\ref{re: 2.2}). From the proof of Lemma \ref{blow-up:le:5.1}, we
know that it is important to handle the second term of the right
hand side of (\ref{blow-up:5.2}) where $\mathrm{div}u$ and
$\nabla|u|$ are involved. On the other hand, the second term of the
left hand side of (\ref{blow-up:5.2}), where $|\nabla u|^2$,
$|\mathrm{div}u|^2$ and $\big|\nabla|u|\big|^2$ are involved, is not
enough to absorb the second term of the right under the physical
restrictions of the viscosities. For the term $|\nabla u|^2$ on the
left of (\ref{blow-up:5.2}), it is natural to get $|\nabla
u|^2\ge\big|\nabla |u|\big|^2$, which makes some additional good
information on $|\nabla u|^2$ lose cf. \cite{Huang-Li-Xin: Serrin,
Sun-Wang-Zhang}. The crucial ingredient to relax the additional
restrictions to $\frac{29\mu}{3}>\lambda$ is that we observe
$$
|\nabla
u|^2=|u|^2\left|\nabla\left(\frac{u}{|u|}\right)\right|^2+\big|\nabla|u|\big|^2,
$$ for $|u|>0$,
and thus
$$
\int_{{\mathbb{R}^3}\cap\{|u|>0\}}|u|^{r-2}|\nabla
u|^2\geq\left(1+\phi(\varepsilon_1,r)\right)\int_{{\mathbb{R}^3}\cap\{|u|>0\}}|u|^{r-2}\big|\nabla|u|\big|^2,
$$
if
$$\label{blow-up:5.8}\int_{{\mathbb{R}^3}\cap\{|u|>0\}}|u|^r\left|\nabla\left(\frac{u}{|u|}\right)\right|^2\geq
\phi(\varepsilon_1,r)\int_{{\mathbb{R}^3}\cap\{|u|>0\}}|u|^{r-2}\big|\nabla|u|\big|^2,$$ for some positive function $\phi(\varepsilon_1,r)$ near $r=3$. For more details, please see Lemma \ref{blow-up:le:5.1}.

(2) In \cite{Choe-Bum}, the authors obtain the upper bound and the positive lower bound of the
density in $\mathbb{T}^3$ under the assumptions $\inf\rho_0>0$, $\mu+\lambda=0$ and $\|\rho\|_{L^\infty(0,T;L^{q_0})}$ is bounded, for some $q_0>0$ large enough and for some $T>0$. From the physical points of view, $\mu+\lambda>0$ seems more natural, since we know that $\mu>0$ and $2\mu+3\lambda\ge0$ deduce $\mu+\lambda>0$. In Theorem \ref{th:1.1}, we only assume $\mu>0$, $2\mu+3\lambda\ge0$ and $\frac{29\mu}{3}>\lambda$.

(3) By (\ref{N-S})$_1$, we known $\|\rho\|_{L^1}=\|\rho_0\|_{L^1}$. It follows from the standard interpolation inequality that the bound of $\|\rho\|_{L^\infty(0,T; L^\infty)}$ yields that $\|\rho\|_{L^\infty(0,T; L^{q_2})}$ is bounded for any $q_2\in(1,\infty)$.  Thus, the blow-up criterion (\ref{result}) is an extension towards (\ref{Huang-Xin and Sun-Wang-Zhang criterion}) in \cite{Huang-Li-Xin: Serrin, Sun-Wang-Zhang}.

\subsection {Full compressible N-S: a blow-up criterion in terms of the upper bound of the density and the temperature}
\setcounter{equation}{0} \setcounter{theorem}{0}
The constants $a$, $C_0$ and $\kappa$ in the equations play no roles in the analysis, we assume $a=C_0=\kappa=1$ henceforth. If the solutions are regular enough (such as strong solutions), (\ref{full N-S}) is equivalence to the following system which is very usefull in the proofs of the main theorems:
 \be\label{full N-S+1}
\begin{cases}
\rho_t+\nabla \cdot (\rho u)=0, \\
\rho  u_t+\rho u\cdot\nabla u+\nabla
P(\rho,\theta)=\mu\Delta u+(\mu+\lambda)\nabla\mathrm{div}u,\\
\rho \theta_t+\rho u\cdot\nabla\theta+\rho\theta\mathrm{div}u=\frac{\mu}{2}\left|\nabla u+(\nabla u)^\prime\right|^2+\lambda(\mathrm{div}u)^2+\Delta\theta,\ \mathrm{in}\ \mathbb{R}^3.
\end{cases}
\ee
System (\ref{full N-S+1}) is supplemented with initial conditions
\be\label{non-initial}
(\rho, u, \theta)|_{t=0}=(\rho_0, u_0, \theta_0),\ x\in\mathbb{R}^3,
\ee
with
\be\label{non-boundary}
\rho(x,t)\rightarrow0,\ u(x,t)\rightarrow0,\ \theta(x,t)\rightarrow0,\ \mathrm{as}\ |x|\rightarrow\infty,\ \mathrm{for}\ t\ge0.
\ee
We give the definition of strong solutions to (\ref{full N-S+1}) throughout the paper.
\begin{definition}(Strong solution) For $T>0$, $(\rho, u, \theta)$ is called a strong solution to the compressible Navier-Stokes equations (\ref{full N-S+1})-(\ref{non-boundary}) in $\mathbb{R}^3\times [0,T]$, if for some $q\in (3, 6]$,
\bex &0\le \rho\in
C([0,T];W^{1,q}\cap H^1\cap L^1),\ \rho_t\in C([0,T];L^2\cap L^q),&\\
& (u, \theta)\in C([0,T];D^2\cap D^1_0)\cap L^2(0,T;D^{2,q}),\ (u_t, \theta_t)\in
L^2(0,T;D^1_0),&\\&
 (\sqrt{\rho} u_t, \sqrt{\rho} \theta_t)\in
L^\infty(0,T;L^2),& \eex
and $(\rho,u,\theta)$ satisfies (\ref{full N-S+1}) a.e. in $\mathbb{R}^3\times (0,T]$.
\end{definition}
Our main result for the full compressible Navier-Stokes equations is stated as follows:
\begin{theorem}\label{non-th:1.1}
Assume $\rho_0\geq0$,
 $\rho_0\in H^1\cap W^{1,q}\cap L^1$, for some $q\in(3,6]$, $(u_0, \theta_0)\in D^2\cap D_0^1$, and the
following compatibility conditions are satisfied:
\beq\label{non-compatibility}
\begin{cases}
\mu\Delta u_0+(\mu+\lambda)\nabla\mathrm{div} u_0-\nabla P(\rho_0, \theta_0)=\sqrt{\rho}_0g_1,\\
\kappa\Delta\theta_0+\frac{\mu}{2}\left|\nabla u_0+(\nabla u_0)^\prime\right|^2+\lambda(\mathrm{div}u_0)^2=\sqrt{\rho_0}g_2,\ x\in\mathbb{R}^3,
\end{cases}
\eeq for some $g_i\in L^2$, $i=1,2$. Let $(\rho, u, \theta)$ be a strong solution to (\ref{full N-S+1})-(\ref{non-boundary}) in $\mathbb{R}^3\times[0, T]$. If $0<T^*<+\infty$ is the maximum time of existence of the strong solution,
then
\be\label{non-result}
\lim\sup\limits_{T\nearrow T^*}\left(\|\rho\|_{L^\infty(0,T; L^\infty)}+\|\theta\|_{L^\infty(0,T; L^\infty)}\right)=\infty,
\ee provided $3\mu>\lambda$.
\end{theorem}
\begin{remark}  Under the conditions of Theorem \ref{non-th:1.1},
the local existence of the strong solutions was obtained in \cite{Cho-Kim; perfect gas}. Thus, the assumption $T^*>0$ makes sense.
\end{remark}
\begin{remark}  Theorem \ref{non-th:1.1} is also valid for more general pressure law, such as $P=a\rho\theta+a_1\rho^\gamma$. Whether the similar result as in Theorem \ref{th:1.1} could be obtained for the full compressible Navier-Stokes equations is still unknown.
\end{remark}

 Before Theorem \ref{non-th:1.1}, there are several results on the blow-up criterions of strong solutions to (\ref{full N-S+1}), please refer for instance to \cite{Fan-Jiang-Ou, Fang-Zi-Zhang, Sun-Wang-Zhang 1} and references therein for initial boundary value problems. In particular,

$\bullet$ Fan-Jiang-Ou (\cite{Fan-Jiang-Ou}, 3D)
\be\label{Fan-Jiang-Ou criterion}
\lim\sup\limits_{t\nearrow
T^\star}\left(\|\theta\|_{L^\infty(0,t;L^\infty)}+\|\nabla
u\|_{L^1(0,t;L^\infty)}\right)=\infty,
\ee
provided $7\mu>\lambda$. Here the appearance of vacuum is allowed.

It is well-known that the bound of $\|\nabla
u\|_{L^1(0,t;L^\infty)}$ yields that $\|\rho\|_{L^\infty(0,t;L^\infty)}$ is bounded (see (2.2) in \cite{Fan-Jiang-Ou}), if the initial density is bounded. When $\|\nabla
u\|_{L^1(0,t;L^\infty)}$ in (\ref{Fan-Jiang-Ou criterion}) is replaced by the upper bound of the density, the following blow-up criterions were obtained:

$\bullet$ Fang-Zi-Zhang (\cite{Fang-Zi-Zhang}, 2D)
\be\label{Fang-Zi-Zhang criterion}
\lim\sup\limits_{t\nearrow
T^\star}\left(\|\theta\|_{L^\infty(0,t;L^\infty)}+\|\rho\|_{L^\infty(0,t;L^\infty)}\right)=\infty,
\ee where the appearance of vacuum is allowed;

$\bullet$ Sun-Wang-Zhang (\cite{Sun-Wang-Zhang 1}, 3D)
\be\label{Sun-Wang-Zhang criterion}
\lim\sup\limits_{t\nearrow
T^\star}\left(\|\theta\|_{L^\infty(0,t;L^\infty)}+\|\rho\|_{L^\infty(0,t;L^\infty)}+\left\|\frac{1}{\rho}\right\|_{L^\infty(0,t;L^\infty)}\right)=\infty,
\ee provided $7\mu>\lambda$.

We would like to point out that an analogous blow-up criterion of (\ref{non-result}) for
the isentropic compressible Naiver-Stokes equation (i.e. (\ref{Huang-Xin and Sun-Wang-Zhang criterion})) in $\mathbb{R}^3$, under the assumption $7\mu>\lambda$, has been previously established by Huang-Li-Xin \cite{Huang-Li-Xin: Serrin} and Sun-Wang-Zhang \cite{Sun-Wang-Zhang}. In \cite{Huang-Li-Xin: Serrin,Sun-Wang-Zhang}, the restriction $7\mu>\lambda$ was needed only for the estimate of $\int\rho|u|^{3+\delta}$ where $\delta>0$ is sufficiently small.

We introduce the main ideas of the proof of Theorem \ref{non-th:1.1}.

(1) To get the upper bound of $\int_{\mathbb{R}^3}\rho |u|^r$, we apply the ideas of the proof of Lemma \ref{blow-up:le:5.1} so that we can get a restriction of $\mu$ and $\lambda$ as better as possible. As a byproduct, we also get the upper bound of $\int_0^t\int_{\mathbb{R}^3}|u|^{r-2}|\nabla u|^2$, which is very crucial in the proof of $L^\infty_tL^2_x$ of $\nabla u$ (see Lemma \ref{non-le: int nabla u}). Here we take $r=4$ because we have to deal with the difficulties caused by the strong nonlinearities in the temperature equation, such as the terms $\frac{\mu}{2}\left|\nabla u+(\nabla u)^\prime\right|^2$ and $\lambda(\mathrm{div}u)^2$ in (\ref{full N-S+1})$_3$, which leads to the restriction $3\mu>\lambda$.

(2) As it was pointed out in \cite{Sun-Wang-Zhang 1} that to deal with the essential difficulties due to the highly nonlinear terms $\left|\nabla u+(\nabla u)^\prime\right|^2$ and $|\mathrm{div}u|^2$ in the temperature equation, Sun-Wang-Zhang used the ideas of Hoff \cite{Hoff ARMA} to get the upper bounds of $L^\infty_t H^s_x$ of $u$ for $s\in(0,1)$, which requires the upper bound of $\frac{1}{\rho}$. Here we do not require the upper bound of $\frac{1}{\rho}$ so that the appearance of vacuum is allowed, because we use the fact $P_t=(\rho E)_t-\left(\frac{\rho|u|^2}{2}\right)_t$, (\ref{full N-S})$_3$ and integration by parts such that
\bex\begin{split}
-\int_{\mathbb{R}^3}P_tG=&-\int_{\mathbb{R}^3}(\rho E)_tG+\cdots=-\int_{\mathbb{R}^3}\mathrm{div}\Big[\Big(\mu\left(\nabla  u+(\nabla u)^\prime\right)+\lambda
\di  u I_N\Big)u\Big]G+\cdots\\ =&\int_{\mathbb{R}^3}\Big[\Big(\mu\left(\nabla  u+(\nabla u)^\prime\right)+\lambda
\di  u I_N\Big)u\Big]\cdot\nabla G+\cdots\\ \le& C\big\||u||\nabla u|\big\|_{L^2}\|\nabla G\|_{L^2}+\cdots,
\end{split}
\eex where $G=(2\mu+\lambda)\mathrm{div}u-P$ is the effective viscous flux which plays an important role in the proofs. For more details, please see (\ref{I 3})-(\ref{I 3+1}) in the proof of Lemma \ref{non-le: int nabla u}.

(3) The nonliear terms $\left|\nabla u+(\nabla u)^\prime\right|^2$ and $|\mathrm{div}u|^2$ in (\ref{full N-S+1})$_3$ could be handled for two space dimension when the blow-up criterion (\ref{Fang-Zi-Zhang criterion}) was established with vacuum, because 2-D Gagliardo-Nirenberg inequality has better properties than 3-D. See \cite{Fang-Zi-Zhang} for more details.
\section{Proof of Theorem \ref{th:1.1}}
\setcounter{equation}{0} \setcounter{theorem}{0}
\renewcommand{\theequation}{\thesection.\arabic{equation}}
\renewcommand{\thetheorem}{\thesection.\arabic{theorem}}
Let $0<T^*<\infty$ be the maximum time of existence of strong solution $(\rho, u)$
to (\ref{N-S})-(\ref{boundary}). Namely,  $(\rho, u)$ is a strong
solution to (\ref{N-S})-(\ref{boundary}) in $\mathbb{R}^3\times [0, T]$ for any $0<T<T^*$, but not
a strong solution in $\mathbb{R}^3\times [0, T^*]$.  Suppose that (\ref{result}) were false,
i.e.
\beq\label{2.1}
M:=\|\rho\|_{L^\infty(0,T^*; L^{q_1})}<\infty. \eeq
The goal is to show that under the assumption (\ref{2.1}),
there is a bound $C>0$ depending only on $M, \rho_0, u_0,\mu,\lambda$, and $T^*$ such that
\beq\label{uniform_est1}
\sup_{0\le t<T^*}\left[\max_{l=2, q}(\|\rho\|_{W^{1,l}}+\|\rho_t\|_{L^l})
+\|\sqrt{\rho}u_t\|_{L^2}+\|\nabla u\|_{H^1}\right]\le C,
\eeq
and
\beq\label{uniform_est2}
\int_{0}^{T^*}\left(\|u_t\|_{D^1}^2+\|u\|_{D^{2,q}}^2\right)\,dt\le C.
\eeq
With (\ref{uniform_est1}) and (\ref{uniform_est2}), it is easy to show without much difficulties
that $T^*$ is not the maximum time, which is the desired contradiction.

Throughout the rest of the section, we denote by $C$ a generic constant depending only on
$\rho_0$, $u_0$, $T^*$, $M$, $\lambda$, $\mu$. We denote by
$$ A\lesssim B $$ if there exists a generic constant $C$ such that $A\leq C B$.

\begin{lemma}\label{blow-up:le:5.1}
Under the conditions of Theorem \ref{th:1.1} and (\ref{2.1}), if $\frac{29\mu}{3}>\lambda$, there
exists $r\in(3,\frac{7}{2})$ such that \bex \sup\limits_{0\leq t\leq
T}\int_{\mathbb{R}^3} \rho|u|^r dx\leq C, \eex for any $T\in[0,T^*)$.
\end{lemma}
\begin{remark}\label{re: 2.2}
Lemma \ref{blow-up:le:5.1} is also true for bounded domains. This lemma relaxes the restriction $7\mu>\lambda$ in \cite{Huang-Li-Xin: Serrin, Sun-Wang-Zhang} to $\frac{29\mu}{3}>\lambda$.  It is easy to verify that Lemma \ref{blow-up:le:5.1} is also true if $P=R\rho\theta$ for a constant $R>0$ and $\theta$ is bounded. Thus, as a byproduct of the paper, the restriction $7\mu>\lambda$ in \cite{Fan-Jiang-Ou, Sun-Wang-Zhang 1} could be relaxed to $\frac{29\mu}{3}>\lambda$ for the full compressible Navier-Stokes equations. In this sense, this lemma extends the results in \cite{Fan-Jiang-Ou, Huang-Li-Xin: Serrin, Sun-Wang-Zhang, Sun-Wang-Zhang 1}.
\end{remark}
\pf  Multiplying (\ref{N-S})$_2$ by
$r|u|^{r-2}u$, and integrating by parts over $\mathbb{R}^3$, we have
\bex\begin{split}
&\frac{d}{dt}\int_{\mathbb{R}^3} \rho|u|^r+\int_{\mathbb{R}^3} r|u|^{r-2}\left(\mu|\nabla u|^2+(\lambda+\mu)|\mathrm{div}u|^2+\mu(r-2)|\nabla|u||^2\right)\\
=&r\int_{\mathbb{R}^3}
\mathrm{div}(|u|^{r-2}u)P-r(r-2)(\mu+\lambda)\int_{\mathbb{R}^3}\mathrm{div}u|u|^{r-3}u\cdot\nabla|u|.\end{split}
\eex Thus,
\beq\label{blow-up:5.2}\begin{split}
&\frac{d}{dt}\int_{\mathbb{R}^3}
\rho|u|^r+\int_{{\mathbb{R}^3}\cap\{|u|>0\}} r|u|^{r-2}\left(\mu|\nabla
u|^2+(\lambda+\mu)|\mathrm{div}u|^2+\mu(r-2)|\nabla|u||^2\right)\\
=&r\int_{{\mathbb{R}^3}\cap\{|u|>0\}}
\mathrm{div}(|u|^{r-2}u)P-r(r-2)(\mu+\lambda)\int_{{\mathbb{R}^3}\cap\{|u|>0\}}\mathrm{div}u|u|^{r-3}u\cdot\nabla|u|.
\end{split}\eeq
For any given $\varepsilon_1\in (0, 1)$, we define a nonnegative function which will be decided in {\bf Case 2} as follows:
$$
\phi(\varepsilon_1,r)=\left\{\begin{array}{l}
\frac{\mu\varepsilon_1 (r-1)}{3\left(-\frac{4\mu}{3}-\lambda+
\frac{r^2(\mu+\lambda)}{4(r-1)}\right)}, \ \ \ {\rm if} \ \ \ \frac{r^2(\mu+\lambda)}{4(r-1)}-\frac{4\mu}{3}-\lambda>0, \\
[3mm] 0, \ \ \ {\rm otherwise}.
\end{array}
\right.
$$
{\bf Case 1:}
\be\label{blow-up:5.8}\int_{{\mathbb{R}^3}\cap\{|u|>0\}}|u|^r\left|\nabla\left(\frac{u}{|u|}\right)\right|^2>\phi(\varepsilon_1,r)\int_{{\mathbb{R}^3}\cap\{|u|>0\}}|u|^{r-2}\big|\nabla|u|\big|^2.\ee
 A direct calculation gives for $|u|>0$\be\label{blow-up:5.4}
|\nabla
u|^2=|u|^2\left|\nabla\left(\frac{u}{|u|}\right)\right|^2+\big|\nabla|u|\big|^2,
\ee which plays a important role in the proof.

By (\ref{blow-up:5.2}), we have \bex
\begin{split}&\frac{d}{dt}\int_{\mathbb{R}^3} \rho|u|^r+\int_{{\mathbb{R}^3}\cap\{|u|>0\}} r|u|^{r-2}\left(\mu|\nabla u|^2+(\lambda+\mu)|\mathrm{div}u|^2+\mu(r-2)\big|\nabla|u|\big|^2\right)\\
=&r\int_{{\mathbb{R}^3}\cap\{|u|>0\}}
\mathrm{div}(|u|^{r-2}u)P-r(r-2)(\mu+\lambda)\int_{{\mathbb{R}^3}\cap\{|u|>0\}}\mathrm{div}u|u|^\frac{r-2}{2}|u|^\frac{r-4}{2}u\cdot\nabla|u|\\
\le&C\int_{{\mathbb{R}^3}\cap\{|u|>0\}}\rho^{\gamma-\frac{r-2}{2r}}
\rho^{\frac{r-2}{2r}} |u|^{r-2}|\nabla
u|+r(\mu+\lambda)\int_{{\mathbb{R}^3}\cap\{|u|>0\}}|u|^{r-2}|\mathrm{div}u|^2\\&+\frac{r(r-2)^2(\mu+\lambda)}{4}\int_{{\mathbb{R}^3}\cap\{|u|>0\}}|u|^{r-2}\big|\nabla|u|\big|^2,
\end{split}\eex where we have used Cauchy inequality.
Thus, \beq\label{blow-up:5.9}
\begin{split}&\frac{d}{dt}\int_{\mathbb{R}^3}
\rho|u|^r+\int_{{\mathbb{R}^3}\cap\{|u|>0\}} \mu r|u|^{r-2}|\nabla
u|^2+\mu(r-2)r\int_{{\mathbb{R}^3}\cap\{|u|>0\}}|u|^{r-2}\big|\nabla|u|\big|^2
\\ \le&C\int_{{\mathbb{R}^3}\cap\{|u|>0\}}
\rho^{\gamma-\frac{r-2}{2r}}\rho^{\frac{r-2}{2r}} |u|^{r-2}|\nabla
u|+\frac{r(r-2)^2(\mu+\lambda)}{4}\int_{{\mathbb{R}^3}\cap\{|u|>0\}}|u|^{r-2}\big|\nabla|u|\big|^2.
\end{split}\eeq
By (\ref{blow-up:5.4}), (\ref{blow-up:5.9}), Cauchy
inequality, and H\"older inequality, for any $\varepsilon_0\in(0,1)$, we have \bex
\begin{split}&\frac{d}{dt}\int_{\mathbb{R}^3}
\rho|u|^r+\int_{{\mathbb{R}^3}\cap\{|u|>0\}} \mu r|u|^{r-2}\big|\nabla
|u|\big|^2+\int_{{\mathbb{R}^3}\cap\{|u|>0\}} \mu r|u|^r\left|\nabla
\left(\frac{u}{|u|}\right)\right|^2
\\&+\mu(r-2)r\int_{{\mathbb{R}^3}\cap\{|u|>0\}}|u|^{r-2}\big|\nabla|u|\big|^2
\\ \le& C\int_{{\mathbb{R}^3}\cap\{|u|>0\}}
\rho^{\gamma-\frac{r-2}{2r}}\rho^{\frac{r-2}{2r}} |u|^{r-2}\big|\nabla
|u|\big|+C\int_{{\mathbb{R}^3}\cap\{|u|>0\}}\rho^{\gamma-\frac{r-2}{2r}}
\rho^{\frac{r-2}{2r}} |u|^{r-1}\left|\nabla
\left(\frac{u}{|u|}\right)\right|\\&+\frac{r(r-2)^2(\mu+\lambda)}{4}\int_{{\mathbb{R}^3}\cap\{|u|>0\}}|u|^{r-2}\big|\nabla|u|\big|^2\\
\le& C\int_{{\mathbb{R}^3}\cap\{|u|>0\}}
\rho^{\gamma-\frac{r-2}{2r}}\rho^{\frac{r-2}{2r}} |u|^{r-2}\big|\nabla
|u|\big|+\mu r\varepsilon_0\int_{{\mathbb{R}^3}\cap\{|u|>0\}}
|u|^r\left|\nabla \left(\frac{u}{|u|}\right)\right|^2\\ &+\frac{C}{4\mu
r\varepsilon_0}\left(\int_{\mathbb{R}^3} \rho |u|^r\right)^\frac{r-2}{r}\left(\int_{\mathbb{R}^3}
\rho^{\frac{(2\gamma-1)r}{2}+1}\right)^\frac{2}{r}
+\frac{r(r-2)^2(\mu+\lambda)}{4}\int_{{\mathbb{R}^3}\cap\{|u|>0\}}|u|^{r-2}\big|\nabla|u|\big|^2.
\end{split}\eex
Combining (\ref{2.1}) and (\ref{blow-up:5.8}), we have \beq\label{blow-up:5.10}
\begin{split}&\frac{d}{dt}\int_{\mathbb{R}^3} \rho|u|^r+r\left[\mu
(1-\varepsilon_0)\phi(\varepsilon_1,r)+\mu
(r-1)-\frac{(r-2)^2(\mu+\lambda)}{4}\right]\int_{{\mathbb{R}^3}\cap\{|u|>0\}}
|u|^{r-2}\big|\nabla |u|\big|^2
\\ \le&C\int_{{\mathbb{R}^3}\cap\{|u|>0\}} \rho^{\gamma-\frac{r-2}{2r}}\rho^{\frac{r-2}{2r}}
|u|^{r-2}\big|\nabla |u|\big|+\frac{C}{4\mu
r\varepsilon_0}\left(\int_{\mathbb{R}^3} \rho
|u|^r\right)^\frac{r-2}{r}.\end{split} \eeq
 ({\bf Sub-Case
1$_1$}):
 If $3\in\{r|\frac{r^2(\mu+\lambda)}{4(r-1)}-\frac{4\mu}{3}-\lambda>0\}$, i.e.,
 $5\mu<3\lambda$, it is easy to get $ [3,\infty)\subset\{r|\frac{r^2(\mu+\lambda)}{4(r-1)}-\frac{4\mu}{3}-\lambda>0\}$.
Therefore, we have
 \be\label{blow-up:5.11}\phi(\varepsilon_1,r)=
\frac{\mu\varepsilon_1 (r-1)}{3\left(-\frac{4\mu}{3}-\lambda+
\frac{r^2(\mu+\lambda)}{4(r-1)}\right)},\ee
 for any $r\in [3,\infty)$.

Denote \be\label{blow-up:5.12}
f(\varepsilon_0,\varepsilon_1,r)=\mu
(1-\varepsilon_0)\phi(\varepsilon_1,r)+\mu
(r-1)-\frac{(r-2)^2(\mu+\lambda)}{4}.
 \ee
Substituting (\ref{blow-up:5.11}) into (\ref{blow-up:5.12}), for
$r\in[3,\infty)$, we have \be\label{blow-up:5.13}
f(\varepsilon_0,\varepsilon_1,r)=\frac{\mu^2\varepsilon_1(1-\varepsilon_0)
(r-1)}{3\left(-\frac{4\mu}{3}-\lambda+
\frac{r^2(\mu+\lambda)}{4(r-1)}\right)}+\mu
(r-1)-\frac{(r-2)^2(\mu+\lambda)}{4}.
 \ee For $(\varepsilon_0,\varepsilon_1,r)=(0,1,3)$, we have
\bex
f(0,1,3)=\frac{16\mu^2}{3\lambda-5\mu}+\frac{7\mu-\lambda}{4}>0,
\eex where we have used $\frac{5\mu}{3}<\lambda<\frac{29}{3}\mu$.

Since $ f(\varepsilon_0,\varepsilon_1,r)$ is continuous
 w.r.t. $(\varepsilon_0, \varepsilon_1, r)$ over $[0, 1]\times [0, 1]\times
[3,\infty)$, there exist $\varepsilon_0, \varepsilon_1 \in (0, 1)$
and $r\in (3,\frac{7}{2})$, such that
$$
f(\varepsilon_0,\varepsilon_1,r)>0.
$$
By (\ref{blow-up:5.10}), Cauchy inequality and H\"older
inequality, we have \bex\begin{split}&\frac{d}{dt}\int_{\mathbb{R}^3}
\rho|u|^{r}+rf(\varepsilon_0,\varepsilon_1,r)\int_{{\mathbb{R}^3}\cap\{|u|>0\}}
|u|^{r-2}\big|\nabla |u|\big|^2
\\ \le&rf(\varepsilon_0,\varepsilon_1,r)\int_{{\mathbb{R}^3}\cap\{|u|>0\}} |u|^{r-2}\big|\nabla
|u|\big|^2+\frac{C}{4rf(\varepsilon_0,\varepsilon_1,r)}\left(\int_{\mathbb{R}^3}
\rho |u|^{r}\right)^\frac{r-2}{r}\left(\int_{\mathbb{R}^3}
\rho^{\frac{(2\gamma-1)r}{2}+1}\right)^\frac{2}{r}\\&+\frac{C}{4\mu
r\varepsilon_0}\left(\int_{\mathbb{R}^3} \rho
|u|^{r}\right)^\frac{r-2}{r}. \end{split}\eex This together with (\ref{2.1}) gives
\beq\label{blow-up:5.14}\frac{d}{dt}\int_{\mathbb{R}^3} \rho|u|^{r} \le
C\left[\frac{1}{f(\varepsilon_0,\varepsilon_1,r)}+\frac{1}{\mu
\varepsilon_0}\right]\left(\int_{\mathbb{R}^3} \rho
|u|^{r}\right)^\frac{r-2}{r}.\eeq
 ({\bf Sub-Case 1$_2$}): if $3\not\in\{r|\frac{r^2(\mu+\lambda)}{4(r-1)}-\frac{4\mu}{3}-\lambda>0\}$, i.e.,
 $5\mu\ge3\lambda$.

In this case, for $r\in(3, \frac{7}{2})$,
it is easy to get \beq\label{blow-up:5.15}
\begin{split}&r\left[\mu
(1-\varepsilon_0)\phi(\varepsilon_1,r)+\mu
(r-1)-\frac{(r-2)^2(\mu+\lambda)}{4}\right]\\&>
3\left(2\mu-\frac{9(\mu+\lambda)}{16}\right)=3\left(\frac{23\mu}{16}-\frac{9\lambda}{16}\right)\\&\ge
3\left(\frac{23\mu}{16}-\frac{15\mu}{16}\right)=\frac{3\mu}{2}.\end{split}
 \eeq
By (\ref{blow-up:5.10}), (\ref{blow-up:5.15}), Cauchy inequality
and H\"older inequality, we have
 \bex
\begin{split}&\frac{d}{dt}\int_{\mathbb{R}^3}
\rho|u|^r+\frac{3\mu}{2}\int_{{\mathbb{R}^3}\cap\{|u|>0\}}
|u|^{r-2}\big|\nabla |u|\big|^2 \\ \le&
C\int_{{\mathbb{R}^3}\cap\{|u|>0\}}\rho^{\gamma-\frac{r-2}{2r}}
\rho^{\frac{r-2}{2r}} |u|^{r-2}\big|\nabla |u|\big|+\frac{C}{4\mu
r\varepsilon_0}\left(\int_{\mathbb{R}^3} \rho |u|^r\right)^\frac{r-2}{r}\\
\le&\frac{3\mu}{2}\int_{{\mathbb{R}^3}\cap\{|u|>0\}}
|u|^{r-2}\big|\nabla |u|\big|^2 +C\left(\int_{\mathbb{R}^3} \rho
|u|^r\right)^\frac{r-2}{r}\left(\int_{\mathbb{R}^3}
\rho^{\frac{(2\gamma-1)r}{2}+1}\right)^\frac{2}{r}+\frac{C}{4\mu
r\varepsilon_0}\left(\int_{\mathbb{R}^3} \rho |u|^r\right)^\frac{r-2}{r}. \end{split}\eex Therefore,
\be\label{blow-up:5.16} \frac{d}{dt}\int_{\mathbb{R}^3} \rho|u|^r\le
C\left(\int_{\mathbb{R}^3} \rho |u|^r\right)^\frac{r-2}{r}, \ee where we have used (\ref{2.1}).

By (\ref{blow-up:5.14}) and (\ref{blow-up:5.16}), for {\bf Case
1}, we conclude that if $\lambda<\frac{29}{3}\mu$ and
(\ref{blow-up:5.8}) are satisfied, the
following estimate can be obtained \be\label{blow-up:5.17}
\frac{d}{dt}\int_{\mathbb{R}^3}\rho|u|^{r}\le C\left(\int_{\mathbb{R}^3} \rho
|u|^{r}\right)^\frac{r-2}{r},\ee for some constants $C>0$ and $r\in(3, \frac{7}{2})$.\\

{\bf Case 2:} if
\be\label{blow-up:5.3}\int_{{\mathbb{R}^3}\cap\{|u|>0\}}|u|^r\left|\nabla\left(\frac{u}{|u|}\right)\right|^2\le\phi(\varepsilon_1,r)\int_{{\mathbb{R}^3}\cap\{|u|>0\}}|u|^{r-2}\big|\nabla|u|\big|^2.\ee
 A direct calculation gives for $|u|>0$ \be\label{blow-up:5.5}
\mathrm{div}u=|u|\mathrm{div}\left(\frac{u}{|u|}\right)+\frac{u\cdot\nabla|u|}{|u|}.
\ee
By (\ref{blow-up:5.2}) and (\ref{blow-up:5.5}), we have \bex
\begin{split}&\frac{d}{dt}\int_{\mathbb{R}^3} \rho|u|^r+\int_{{\mathbb{R}^3}\cap\{|u|>0\}} r|u|^{r-2}\left(\mu|\nabla u|^2+(\lambda+\mu)|\mathrm{div}u|^2+\mu(r-2)\big|\nabla|u|\big|^2\right)\\
=& r\int_{{\mathbb{R}^3}\cap\{|u|>0\}}
\mathrm{div}(|u|^{r-2}u)P-r(r-2)(\mu+\lambda)\int_{{\mathbb{R}^3}\cap\{|u|>0\}}|u|^{r-2}u\cdot\nabla|u|\mathrm{div}\left(\frac{u}{|u|}\right)
\\&-r(r-2)(\mu+\lambda)\int_{{\mathbb{R}^3}\cap\{|u|>0\}}|u|^{r-4}\big|u\cdot\nabla|u|\big|^2.
\end{split}\eex
This gives \beq\label{dt rho u r}
\begin{split}&\frac{d}{dt}\int_{\mathbb{R}^3} \rho|u|^r+\int_{{\mathbb{R}^3}\cap\{|u|>0\}} r|u|^{r-4}G
= r\int_{{\mathbb{R}^3}\cap\{|u|>0\}} \mathrm{div}(|u|^{r-2}u)P,
\end{split}\eeq where
\bex\begin{split} G=&\mu|u|^2|\nabla
u|^2+(\lambda+\mu)|u|^2|\mathrm{div}u|^2+\mu(r-2)|u|^2\big|\nabla|u|\big|^2
\\&+(r-2)(\mu+\lambda)|u|^2u\cdot\nabla|u|\mathrm{div}\left(\frac{u}{|u|}\right)+(r-2)(\mu+\lambda)\big|u\cdot\nabla|u|\big|^2.\end{split}
\eex
To let $\int_{{\mathbb{R}^3}\cap\{|u|>0\}} r|u|^{r-4}G$ become a {\it good} term, we shall consider $G$ first.
\bex\begin{split}
G=&\mu|u|^2\left(|u|^2\left|\nabla\left(\frac{u}{|u|}\right)\right|^2+\big|\nabla|u|\big|^2\right)+(\mu+\lambda)|u|^2\left(|u|\mathrm{div}\left(\frac{u}{|u|}\right)
+\frac{u\cdot\nabla
|u|}{|u|}\right)^2\\&+\mu(r-2)|u|^2\big|\nabla
|u|\big|^2+(r-2)(\mu+\lambda)|u|^2u\cdot\nabla|u|\mathrm{div}\left(\frac{u}{|u|}\right)+(r-2)(\mu+\lambda)\big|u\cdot\nabla
|u|\big|^2\\=&\mu|u|^4\left|\nabla\left(\frac{u}{|u|}\right)\right|^2+\mu(r-1)|u|^2\big|\nabla|u|\big|^2+(r-1)(\mu+\lambda)\big|u\cdot\nabla
|u|\big|^2\\&+r(\mu+\lambda)
|u|^2u\cdot\nabla|u|\mathrm{div}\left(\frac{u}{|u|}\right)+(\mu+\lambda)|u|^4\left(\mathrm{div}\left(\frac{u}{|u|}\right)\right)^2\\=&
\mu|u|^4\left|\nabla\left(\frac{u}{|u|}\right)\right|^2+\mu(r-1)|u|^2\big|\nabla|u|\big|^2
+(r-1)(\mu+\lambda)\left(u\cdot\nabla|u|+\frac{r}{2(r-1)}|u|^2\mathrm{div}\left(\frac{u}{|u|}\right)\right)^2
\\&+(\mu+\lambda)|u|^4\left(\mathrm{div}\left(\frac{u}{|u|}\right)\right)^2-
\frac{r^2(\mu+\lambda)}{4(r-1)}|u|^4\left(\mathrm{div}\left(\frac{u}{|u|}\right)\right)^2.
\end{split}
\eex This, combining the fact
$$
\left|\mathrm{div}\left(\frac{u}{|u|}\right)\right|^2\le3\left|\nabla\left(\frac{u}{|u|}\right)\right|^2,
$$ deduces
\bex\begin{split} G\ge&
\mu|u|^4\left|\nabla\left(\frac{u}{|u|}\right)\right|^2+\mu(r-1)|u|^2\big|\nabla|u|\big|^2
+\left(\mu+\lambda-
\frac{r^2(\mu+\lambda)}{4(r-1)}\right)|u|^4\left(\mathrm{div}\left(\frac{u}{|u|}\right)\right)^2\\
\ge&\frac{\mu}{3}|u|^4\left(\mathrm{div}\left(\frac{u}{|u|}\right)\right)^2+\left(\mu+\lambda-
\frac{r^2(\mu+\lambda)}{4(r-1)}\right)|u|^4\left(\mathrm{div}\left(\frac{u}{|u|}\right)\right)^2+\mu(r-1)|u|^2\big|\nabla|u|\big|^2\\
=&\left(\frac{4\mu}{3}+\lambda-
\frac{r^2(\mu+\lambda)}{4(r-1)}\right)|u|^4\left(\mathrm{div}\left(\frac{u}{|u|}\right)\right)^2+\mu(r-1)|u|^2\big|\nabla|u|\big|^2.
\end{split}
\eex
Thus, \bex\begin{split} \int_{{\mathbb{R}^3}\cap\{|u|>0\}}
r|u|^{r-4}G\ge&r\left(\frac{4\mu}{3}+\lambda-
\frac{r^2(\mu+\lambda)}{4(r-1)}\right)\int_{{\mathbb{R}^3}\cap\{|u|>0\}}
|u|^{r}\left(\mathrm{div}\left(\frac{u}{|u|}\right)\right)^2\\&+\mu
r(r-1)\int_{{\mathbb{R}^3}\cap\{|u|>0\}}
|u|^{r-2}\big|\nabla|u|\big|^2\\
\ge&3r\left(\frac{4\mu}{3}+\lambda-
\frac{r^2(\mu+\lambda)}{4(r-1)}\right)\phi(\varepsilon_1,r)\int_{{\mathbb{R}^3}\cap\{|u|>0\}}
|u|^{r-2}\big|\nabla|u|\big|^2\\&+\mu
r(r-1)\int_{{\mathbb{R}^3}\cap\{|u|>0\}}
|u|^{r-2}\big|\nabla|u|\big|^2\\=&\left[3r\left(\frac{4\mu}{3}+\lambda-
\frac{r^2(\mu+\lambda)}{4(r-1)}\right)\phi(\varepsilon_1,r)+\mu
r(r-1)\right]\int_{{\mathbb{R}^3}\cap\{|u|>0\}}
|u|^{r-2}\big|\nabla|u|\big|^2,
\end{split}
\eex
where we have used (\ref{blow-up:5.3}).

Putting all these estimates into (\ref{dt rho u r}), we have
 \bex
\begin{split}&\frac{d}{dt}\int_{\mathbb{R}^3} \rho|u|^r+\left[3r\left(\frac{4\mu}{3}+\lambda-
\frac{r^2(\mu+\lambda)}{4(r-1)}\right)\phi(\varepsilon_1,r)+\mu
r(r-1)\right]\int_{{\mathbb{R}^3}\cap\{|u|>0\}}
|u|^{r-2}\big|\nabla|u|\big|^2\\ \le&
C\int_{{\mathbb{R}^3}\cap\{|u|>0\}}
\rho^{\gamma-\frac{r-2}{2r}}\rho^{\frac{r-2}{2r}}|u|^{r-2}|\nabla u|\\
\le&\varepsilon\int_{{\mathbb{R}^3}\cap\{|u|>0\}}|u|^{r-2}|\nabla
u|^2+\frac{C}{\varepsilon}\left(\int_{{\mathbb{R}^3}\cap\{|u|>0\}}
\rho|u|^r\right)^\frac{r-2}{r}\left(\int_{{\mathbb{R}^3}\cap\{|u|>0\}}
\rho^{\frac{(2\gamma-1)r}{2}+1}\right)^\frac{2}{r}\\
\le&\varepsilon\big(1+\phi(\varepsilon_1,r)\big)\int_{{\mathbb{R}^3}\cap\{|u|>0\}}
|u|^{r-2}\big|\nabla|u|\big|^2+\frac{C}{\varepsilon}\left(\int_{{\mathbb{R}^3}\cap\{|u|>0\}}
\rho|u|^r\right)^\frac{r-2}{r},
\end{split}\eex
where we have used Cauchy inequality, H\"older inequality and (\ref{2.1}).

Taking
$\varepsilon=\big(1+\phi(\varepsilon_1,r)\big)^{-1}\left[3r\left(\frac{4\mu}{3}+\lambda-
\frac{r^2(\mu+\lambda)}{4(r-1)}\right)\phi(\varepsilon_1,r)+\mu
r(r-1)\right]$, we have \beq\label{blow-up:5.7}
\begin{split}\frac{d}{dt}\int_{\mathbb{R}^3} \rho|u|^r
\le
\frac{C\big(1+\phi(\varepsilon_1,r)\big)}{\left[3r\left(\frac{4\mu}{3}+\lambda-
\frac{r^2(\mu+\lambda)}{4(r-1)}\right)\phi(\varepsilon_1,r)+\mu
r(r-1)\right]}\left(\int_{{\mathbb{R}^3}} \rho|u|^r\right)^\frac{r-2}{r},
\end{split}\eeq
for $r\in(3,\frac{7}{2})$.

  By (\ref{blow-up:5.17}) and (\ref{blow-up:5.7}), for {\bf Case
1} and {\bf Case 2}, we conclude that if $\lambda<\frac{29}{3}\mu$, there exist some constants $C>0$ and $r\in(3,\frac{7}{2})$ such that  \be\label{blow-up:5.18}
\frac{d}{dt}\int_{\mathbb{R}^3} \rho|u|^{r}\le C\left(\int_{\mathbb{R}^3} \rho
|u|^{r}\right)^\frac{r-2}{r}. \ee
Since
$\frac{r-2}{r}\in(0,1)$, using Young inequality and
Gronwall inequality over (\ref{blow-up:5.18}), we complete the proof of Lemma
\ref{blow-up:le:5.1}.
\endpf \\

From Remark \ref{re: 2.2} and \cite{Huang-Li-Xin: Serrin, Sun-Wang-Zhang}, in order to get (\ref{uniform_est1}) and (\ref{uniform_est2}), it suffices to get the upper bound of $\sup\limits_{0\le t<T^*}\|\rho(t)\|_{L^\infty}$. To do this, Lemma \ref{le: int nabla u} and Lemma \ref{i-le: int rho u t} are needed.
\begin{lemma}\label{le: int nabla u}Under the conditions of Theorem \ref{th:1.1} and (\ref{2.1}), it holds that for any $T\in[0,T^*)$
$$\sup\limits_{0\le t\le T}\int_{\mathbb{R}^3}|\nabla u|^2+\int_0^T\int_{\mathbb{R}^3}\rho |\dot{u}|^2\le C,$$ where $\dot{u}=u_t+u\cdot\nabla u$ by the definition of the material derivative.
\end{lemma}
\pf
Multiplying (\ref{N-S})$_2$ by $u_t$, and integrating by parts over $\mathbb{R}^3$, we have
\beq\label{i-dt nabla u 2-1}\begin{split}
&\int_{\mathbb{R}^3}\rho |\u|^2+\frac{1}{2}\frac{d}{dt}\int_{\mathbb{R}^3}\left(\mu|\nabla u|^2+(\mu+\lambda)|\mathrm{div}u|^2\right)\\=&\int_{\mathbb{R}^3}\rho u\cdot\nabla u \cdot \u+\frac{d}{dt}\int_{\mathbb{R}^3}P\mathrm{div}u-\int_{\mathbb{R}^3}P_t\mathrm{div}u\\=&\frac{d}{dt}\int_{\mathbb{R}^3}P\mathrm{div}u-\frac{1}{2(2\mu+\lambda)}
\frac{d}{dt}\int_{\mathbb{R}^3}P^2
-\frac{1}{2\mu+\lambda}\int_{\mathbb{R}^3}P_tG+\int_{\mathbb{R}^3}\rho u\cdot\nabla u \cdot \u\\=&\frac{d}{dt}\int_{\mathbb{R}^3}P\mathrm{div}u-\frac{1}{2(2\mu+\lambda)}
\frac{d}{dt}\int_{\mathbb{R}^3}P^2+\frac{1}{2\mu+\lambda}\int_{\mathbb{R}^3}\big(\mathrm{div}(P u)+(\gamma-1)P\mathrm{div}u\big)G\\&+\int_{\mathbb{R}^3}\rho u\cdot\nabla u \cdot \u\\=&\frac{d}{dt}\int_{\mathbb{R}^3}P\mathrm{div}u-\frac{1}{2(2\mu+\lambda)}
\frac{d}{dt}\int_{\mathbb{R}^3}P^2-
\frac{1}{2\mu+\lambda}\int_{\mathbb{R}^3}P u\cdot\nabla G+\frac{\gamma-1}{2\mu+\lambda}\int_{\mathbb{R}^3}P\mathrm{div}u G\\&+\int_{\mathbb{R}^3}\rho u\cdot\nabla u \cdot \u=\sum\limits_{i=1}^5I_i,
\end{split}
\eeq where $G=(2\mu+\lambda)\mathrm{div}u-P$.

For $I_3$, using H\"older inequality, we have
\beq\label{i-I 3}\begin{split}
I_3\les\int_{\mathbb{R}^3}P |u||\nabla G| \les\|\rho^\frac{1}{r}u\|_{L^r}\left\|\rho^{\gamma-\frac{1}{r}}\right\|_{L^\frac{rp_1}{rp_1-p_1-r}}\|\nabla G\|_{L^{p_1}},
\end{split}
\eeq  for some $p_1\in(1,2)$.

Taking $\mathrm{div}$ on both side of (\ref{N-S})$_2$, we have
\be\label{i-equation of G}
\Delta G=\mathrm{div}(\rho \u).
\ee
From the standard elliptic estimates together with (\ref{2.1}), we have
\beq\label{w 1, p-1 of G}
\begin{split}\|\nabla G\|_{L^{p_1}}\les \|\rho \u\|_{L^{p_1}} \les
\|\sqrt{\rho} \u\|_{L^2}\Big\|\sqrt{\rho}\Big\|_{L^\frac{2p_1}{2-p_1}} \les\|\sqrt{\rho} \u\|_{L^2}.
\end{split}
\eeq
By (\ref{i-I 3}), (\ref{w 1, p-1 of G}), Lemma \ref{blow-up:le:5.1} and (\ref{2.1}), we have
\beq\label{i-I 3+1}\begin{split}
I_3\le&C\|\sqrt{\rho} \u\|_{L^2}.
\end{split}
\eeq
For $I_4$, we have
\beq\label{i-I 4}
\begin{split}
I_4\les&\int_{\mathbb{R}^3}P|\mathrm{div}u| |G|\\ \les&\Big\|P\Big\|_{L^\frac{6p_1}{5p_1-6}}\|\mathrm{div}u\|_{L^2}\Big\|G\Big\|_{L^\frac{3p_1}{3-p_1}}\\ \les&
\|\mathrm{div}u\|_{L^2}\|\nabla G\|_{L^{p_1}},
\end{split}
\eeq
where we have used H\"older inequality, Sobolev inequality and (\ref{2.1}).

Substituting (\ref{w 1, p-1 of G}) into (\ref{i-I 4}), we have
\beq\label{i-I 4+1}
\begin{split}
I_4\le C\|\mathrm{div}u\|_{L^2}\|\sqrt{\rho} \u\|_{L^2}.
\end{split}
\eeq
For $I_5$, we have
\beq\label{i-I 5}
\begin{split}
I_5\le& \|\sqrt{\rho} \u\|_{L^2}\|\sqrt{\rho}u\cdot\nabla u\|_{L^2}.
\end{split}
\eeq
Assume $p_2\in (\frac{2r}{r-2},6)$, and let $\frac{3p_2}{3+p_2}<p_1$, we have for any $\varepsilon\in(0,1)$
\beq\label{rho 1/2 u nabla u}\begin{split}
\|\sqrt{\rho} u\cdot\nabla u\|_{L^{2}}\le&\|\rho^\frac{1}{r}u\|_{L^r}\Big\|\rho^{\frac{1}{2}-\frac{1}{r}}\Big\|_{L^\frac{2rp_2}{rp_2-2p_2-2r}}\|\nabla u\|_{L^{p_2}}
\\ \les&\|\mathrm{div} u\|_{L^{p_2}}+\|\mathrm{curl}u\|_{L^{p_2}}\\ \les&\|G\|_{L^{p_2}}+\|\mathrm{curl}u\|_{L^{p_2}}+1\\ \le&\varepsilon\|\nabla G\|_{L^{p_1}}+\varepsilon\|\nabla\mathrm{curl}u\|_{L^{p_1}}+C_\varepsilon\|\nabla u\|_{L^2}+C,
\end{split}
\eeq where we have used H\"older inequality, (\ref{2.1}), Lemma \ref{blow-up:le:5.1}, and the standard interpolation inequality.

Taking $\mathrm{curl}$ on both side of (\ref{N-S})$_2$, we have
\bex
\mu\Delta(\mathrm{curl}u)=\mathrm{curl}(\rho \u).
\eex
Similar to (\ref{w 1, p-1 of G}), we have
\beq\label{w 1, p-1 of curl u}
\begin{split}\|\nabla \mathrm{curl}u\|_{L^{p_1}} \les
\|\sqrt{\rho} \u\|_{L^2}.
\end{split}
\eeq
Substituting (\ref{w 1, p-1 of G}) and (\ref{w 1, p-1 of curl u}) into (\ref{rho 1/2 u nabla u}), we have
\beq\label{rho 1/2 u nabla u+1}\begin{split}
\|\sqrt{\rho} u\cdot\nabla u\|_{L^{2}}\le&\varepsilon C\|\sqrt{\rho} \u\|_{L^2}+C_\varepsilon\|\nabla u\|_{L^2}+C.
\end{split}
\eeq
Substituting (\ref{rho 1/2 u nabla u+1}) into (\ref{i-I 5}), we have
\beq\label{i-I 5+1}
\begin{split}
I_5\le&\varepsilon C\|\sqrt{\rho} \u\|_{L^2}^2+C_\varepsilon\|\nabla u\|_{L^2}^2+C.
\end{split}
\eeq
Putting (\ref{i-I 3+1}), (\ref{i-I 4+1}) and (\ref{i-I 5+1}) into (\ref{i-dt nabla u 2-1}), using Cauchy inequality, and taking $\varepsilon$ sufficiently small, we have
\beq\label{dt nabla u 2-2}\begin{split}
&\frac{1}{2}\int_{\mathbb{R}^3}\rho |\u|^2+\frac{1}{2}\frac{d}{dt}\int_{\mathbb{R}^3}\left(\mu|\nabla u|^2+(\mu+\lambda)|\mathrm{div}u|^2\right)\\ \le&\frac{d}{dt}\int_{\mathbb{R}^3}P\mathrm{div}u-\frac{1}{2(2\mu+\lambda)}
\frac{d}{dt}\int_{\mathbb{R}^3}P^2+C\|\nabla u\|_{L^2}^2+C.
\end{split}
\eeq
Integrating (\ref{dt nabla u 2-2}) over $[0, t]$, and using Cauchy inequality, we have
\bex\begin{split}
&\frac{1}{2}\int_0^t\int_{\mathbb{R}^3}\rho |\u|^2+\frac{1}{2}\int_{\mathbb{R}^3}\left(\mu|\nabla u|^2+(\mu+\lambda)|\mathrm{div}u|^2\right)\\ \le&\int_{\mathbb{R}^3}P\mathrm{div}u+C\int_0^t\|\nabla u\|_{L^2}^2+C\le\frac{\mu+\lambda}{2}\int_{\mathbb{R}^3}|\mathrm{div}u|^2+C\int_{\mathbb{R}^3}P^2+C\int_0^t\|\nabla u\|_{L^2}^2+C.
\end{split}
\eex
This together with Gronwall inequality gives
\bex
\int_0^t\int_{\mathbb{R}^3}\rho |\u|^2+\int_{\mathbb{R}^3}|\nabla u|^2\le C,
\eex for any $t\in[0,T^*)$.
\endpf

\begin{lemma}\label{i-le: int rho u t}Under the conditions of Theorem \ref{th:1.1} and (\ref{2.1}), it holds that for any $T\in[0,T^*)$
$$\sup\limits_{0\le t\le T}\int_{\mathbb{R}^3}\rho |\dot{u}|^2+\int_0^T\int_{\mathbb{R}^3}|\nabla \dot{u}|^2\le C.$$
\end{lemma}
\pf
By the definition of $\dot{u}$, we can write (\ref{N-S})$_2$ as follows:
\beq\label{i-Momentum-1}\rho\dot{u}+\nabla (P(\rho))=\mu\Delta u+(\mu+\lambda)\nabla\mathrm{div}u.\eeq
Differentiating (\ref{i-Momentum-1}) with respect to $t$ and using (\ref{N-S})$_1$, we have
\beq\label{i-rho u t=.}
\begin{split}
\rho\u_t+\rho u\cdot\nabla\u+\nabla P_t
=&\mu\Delta\u+(\mu+\lambda)\nabla\mathrm{div}\u-\mu\Delta(u\cdot\nabla u)-(\mu+\lambda)\nabla\mathrm{div}(u\cdot\nabla u)\\&+\di\Big(\mu\Delta u\otimes u+(\mu+\lambda)\nabla\mathrm{div}u\otimes u-\nabla P\otimes u\Big).
\end{split}
\eeq
Multiplying (\ref{i-rho u t=.}) by $\u$, integrating by parts over $\mathbb{R}^3$, for $t\in(0,T^*)$, we obtain
\beq\label{i-dt rho u t-1}
\begin{split}
&\frac{1}{2}\frac{d}{dt}\int_{\mathbb{R}^3}\rho|\u|^2+\int_{\mathbb{R}^3}\left(\mu|\nabla\u|^2+(\mu+\lambda)|\di\u|^2\right)\\
=&\int_{\mathbb{R}^3}\left( P_t\di\u + u\otimes\nabla P:\nabla\u\right)+\mu\int_{\mathbb{R}^3}\Big(\di(\de u\otimes u)-\de(u\cdot\nabla u)\Big)\cdot\u \\
+&(\mu+\lambda)\int_{\mathbb{R}^3}\Big(\di(\nabla\di u\otimes u)-\nabla\di(u\cdot\nabla u)\Big)\cdot\u
=\sum\ls_{i=1}^3II_i.
\end{split}
\eeq
For $II_1$, using (\ref{N-S})$_1$, we have
\beq\label{i-II 1}
\begin{split}
II_1=&\int_{\mathbb{R}^3}\Big(-\di(Pu)\di\u -(\gamma-1)P\di u\di\u + u\otimes\nabla P:\nabla\u \Big)\\
=&\int_{\mathbb{R}^3}\Big(Pu\cdot\nabla\di\u-(\gamma-1)P\di u\di\u -P(\nabla u)^t:\nabla\u -P u\cdot\nabla\di\u\Big) \\
=&-\int_{\mathbb{R}^3}\Big((\gamma-1)P\di u\di\u + P(\nabla u)^t:\nabla\u \Big)
\les\|P\|_{L^4}\|\nabla u\|_{L^4}\|\nabla \u\|_{L^2}.
\end{split}
\eeq
For $II_2$ and $II_3$, we use the similar arguments as \cite{Huang-Wang-Wen, Huang-Li-Xin: Serrin, Sun-Wang-Zhang, Sun-Wang-Zhang 1}. More precisely, we have
$$\di(\de u\otimes u)-\de(u\cdot\nabla u)
=\nabla_k(\di u\nabla_ku)-\nabla_k(\nabla_ku^j\nabla_ju)-\nabla_j(\nabla_ku^j\nabla_ku).$$
Using integration by parts, we have
\beq\label{i-II 2}
\begin{split}
II_2
=\mu\int\Big(\nabla_k(\di u\nabla_ku)-\nabla_k(\nabla_ku^j\nabla_ju)-\nabla_j(\nabla_ku^j\nabla_ku)\Big)\cdot\u
\les\|\nabla \u\|_{L^2}\|\nabla u\|_{L^4}^2.
\end{split}
\eeq
Similarly, since
$$
\di(\nabla\di u\otimes u)-\nabla\di(u\cdot\nabla u)=\nabla(\nabla_ju^j\nabla_iu^i)-\nabla(\nabla_ju^i\nabla_iu^j)-\nabla_i(\nabla u^i\nabla_ju^j),
$$
we have
\beq\label{i-II 3}
\begin{split}
II_3=
(\mu+\lambda)\int\left(\nabla(\nabla_ju^j\nabla_iu^i)-\nabla(\nabla_ju^i\nabla_iu^j)-\nabla_i(\nabla u^i\nabla_ju^j)\right)\cdot\u
\les\|\nabla \u\|_{L^2}\|\nabla u\|_{L^4}^2.
\end{split}
\eeq
Substituting (\ref{i-II 1}), (\ref{i-II 2}) and (\ref{i-II 3}) into (\ref{i-dt rho u t-1}), and using Cauchy inequality and (\ref{2.1}), we have
\bex
\begin{split}
\frac{1}{2}\frac{d}{dt}\int_{\mathbb{R}^3}\rho|\u|^2+\int_{\mathbb{R}^3}\left(\mu|\nabla\u|^2+(\mu+\lambda)|\di\u|^2\right)
\le\frac{\mu}{2}\|\nabla \u\|_{L^2}^2+C\|\nabla u\|_{L^4}^4+C.
\end{split}
\eex
This gives
\beq\label{i-dt rho u t-2}
\begin{split}
&\frac{d}{dt}\int_{\mathbb{R}^3}\rho|\u|^2+\mu\int_{\mathbb{R}^3}|\nabla\u|^2\\
\le& C\|\nabla u\|_{L^4}^4+C \les\|\mathrm{div} u\|_{L^4}^4+\|\mathrm{curl}u\|_{L^4}^4+1\\ \les&\|G\|_{L^4}^4+\|\mathrm{curl}u\|_{L^4}^4+1\\ \les&
\|G\|_{L^2}^\frac{2(7p_1-12)}{5p_1-6}\|\nabla G\|_{L^{p_1}}^\frac{6p_1}{5p_1-6}+\|\mathrm{curl}u\|_{L^2}^\frac{2(7p_1-12)}{5p_1-6}\|\nabla\mathrm{curl}u\|_{L^{p_1}}^\frac{6p_1}{5p_1-6}+1,
\end{split}
\eeq where we have used Gagliardo-Nirenberg inequality and (\ref{2.1}).

By (\ref{i-dt rho u t-2}), Lemma \ref{le: int nabla u}, (\ref{2.1}), (\ref{w 1, p-1 of G}), (\ref{w 1, p-1 of curl u}) and Young inequality, we have
\beq\label{i-dt rho u t-3}
\begin{split}
\frac{d}{dt}\int_{\mathbb{R}^3}\rho|\u|^2+\mu\int_{\mathbb{R}^3}|\nabla\u|^2 \les
\|\sqrt{\rho}\u\|_{L^2}^\frac{6p_1}{5p_1-6}+1 \les
\|\sqrt{\rho}\u\|_{L^2}^4+1,
\end{split}
\eeq where we have used the fact $\frac{6p_1}{5p_1-6}<4$, since $p_1>\frac{3p_2}{3+p_2}\ge\frac{12}{7}$.

Since $\|\sqrt{\rho}\u\|_{L^2}^2$ is bounded in $L^1(0,T)$ (see Lemma \ref{le: int nabla u}), we apply (\ref{i-dt rho u t-3}) and Gronwall inequality to complete the proof of Lemma \ref{i-le: int rho u t}.
\endpf
\begin{corollary}\label{cor:2.5}Under the conditions of Theorem \ref{th:1.1} and (\ref{2.1}), it holds that for any $T\in[0,T^*)$
$$
\|\nabla G\|_{L^2(0,T; L^\frac{6p_1}{12-5p_1})}\le C.
$$
\end{corollary}
\pf By (\ref{i-equation of G}) and the standard elliptic estimates, together with (\ref{2.1}), H\"older inequality, Sobolev inequality and Lemma \ref{i-le: int rho u t}, we have for any $T\in[0,T^*)$
\bex\begin{split}
\int_0^T\Big\|\nabla G\Big\|_{L^\frac{6p_1}{12-5p_1}}^2\le& C\int_0^T\Big\|\rho\u\Big\|_{L^\frac{6p_1}{12-5p_1}}^2\le C\int_0^T\Big\|\rho\Big\|_{L^\frac{p_1}{2-p_1}}^2\|\u\|_{L^6}^2\le
 C\int_0^T\|\nabla\u\|_{L^2}^2\le C.
\end{split}
\eex
\endpf

\begin{lemma}\label{le: L infty of rho}Under the conditions of Theorem \ref{th:1.1} and (\ref{2.1}), it holds that for any $T\in[0,T^*)$
$$\|\rho\|_{L^\infty(0,T; L^\infty)}\le C.$$
\end{lemma}
\pf
For any $1<p<+\infty$, multiplying (\ref{N-S})$_1$ by $p\rho^{p-1}$
and integrating by parts over $\mathbb{R}^3$, we obtain
\beq\label{rho infty 1}
\begin{split}
\frac{d}{dt}\int_{\mathbb{R}^3}\rho^{p}
=&-\int_{\mathbb{R}^3}\Big(u\cdot\nabla(\rho^{p})+p\rho^{p}\di u\Big)\\
=&(1-p)\int_{\mathbb{R}^3}\rho^{p}\di u
=\frac{1-p}{2\mu+\lambda}\int_{\mathbb{R}^3}\rho^{p}G+\frac{1-p}{2\mu+\lambda}\int_{\mathbb{R}^3}\rho^{p}P\\ \le&\frac{p-1}{2\mu+\lambda}\|G\|_{L^\infty}\int_{\mathbb{R}^3}\rho^{p}.
\end{split}
\eeq
Since $\frac{6p_1}{12-5p_1}>3$, using the standard interpolation inequality, we have
\be\label{G infty}
\|G\|_{L^\infty}\le C\Big\|\nabla G\Big\|_{L^\frac{6p_1}{12-5p_1}}+C\|G\|_{L^2}\le C\Big\|\nabla G\Big\|_{L^\frac{6p_1}{12-5p_1}}+C.
\ee
Substituting (\ref{G infty}) into (\ref{rho infty 1}), we have
\bex\begin{split}\frac{d}{dt}\|\rho\|_{L^p}\le& \frac{C(p-1)}{p}\left(\Big\|\nabla G\Big\|_{L^\frac{6p_1}{12-5p_1}}+1\right)\|\rho\|_{L^p}\\ \le&
C\left(\Big\|\nabla G\Big\|_{L^\frac{6p_1}{12-5p_1}}+1\right)\|\rho\|_{L^p},\end{split}\eex where the constant $C$ is independent of $p$.
This and Corollary \ref{cor:2.5}, together with Gronwall inequality, give
$$\sup_{0\le t\le T}\|\rho(t)\|_{L^p}
\leq \|\rho_0\|_{L^p}\exp\left(C\int_0^{T}\left(\Big\|\nabla G\Big\|_{L^\frac{6p_1}{12-5p_1}}+1\right)dt\right)\le C,$$  for any $T\in[0,T^*)$.
 Let $p$ go to $\infty$, we complete the proof of Lemma \ref{le: L infty of rho}.
\endpf

\section{Proof of Theorem \ref{non-th:1.1}}
\setcounter{equation}{0} \setcounter{theorem}{0}
Let $0<T^*<\infty$ be the maximum time of existence of strong solution $(\rho, u)$
to (\ref{full N-S+1})-(\ref{non-boundary}). Namely,  $(\rho, u)$ is a strong
solution to (\ref{full N-S+1})-(\ref{non-boundary}) in $\mathbb{R}^3\times [0, T]$ for any $0<T<T^*$, but not
a strong solution in $\mathbb{R}^3\times [0, T^*]$. We shall prove Theorem \ref{non-th:1.1} by using a contradiction
argument. Suppose that (\ref{non-result}) were false,
i.e.
\beq\label{non-2.1}
M:=\|\rho\|_{L^\infty(0,T^*; L^\infty)}+\|\theta\|_{L^\infty(0,T^*; L^\infty)}<\infty. \eeq
The goal is to show that under the assumption (\ref{non-2.1}),
there is a bound $C>0$ depending only on $M, \rho_0, u_0, \theta_0, \mu,\lambda, \kappa$, and $T^*$ such that
\beq\label{non-uniform_est1}
\sup_{0\le t<T^*}\left[\max_{l=2, q}(\|\rho\|_{W^{1,l}}+\|\rho_t\|_{L^l})
+\|(\sqrt{\rho}u_t, \sqrt{\rho}\theta_t)\|_{L^2}+\|(\nabla u, \nabla\theta)\|_{H^1}\right]\le C,
\eeq
and
\beq\label{non-uniform_est2}
\int_{0}^{T^*}\left(\|(u_t, \theta_t)\|_{D^1}^2+\|(u, \theta)\|_{D^{2,q}}^2\right)\,dt\le C.
\eeq
With (\ref{non-uniform_est1}) and (\ref{non-uniform_est2}), it is easy to show without much difficulties
that $T^*$ is not the maximum time, which is the desired contradiction.

Throughout the rest of the section, we denote by $C$ a generic constant depending only on
$\rho_0$, $u_0$, $\theta_0$, $T^*$, $M$, $\lambda$, $\mu$, $\kappa$. We denote by
$$A\lesssim B$$ if there exists a generic constant $C$ such that $A\leq C B$.
\begin{lemma}\label{non-le:2.1}
Under the conditions of Theorem \ref{non-th:1.1} and (\ref{non-2.1}), it holds that
\be\label{non-energy inequality}
\begin{cases}
\sup\limits_{0\leq t\leq
T}\int_{\mathbb{R}^3} \rho(|u|^2+\theta^2)\, dx+\int_0^T\int_{\mathbb{R}^3}\left(|\nabla u|^2+|\nabla \theta|^2\right)\,dx\leq C,\\
\sup\limits_{0\leq t\leq
T}\int_{\mathbb{R}^3} \rho\le C,\ \mathrm{for}\ \mathrm{any}\ T\in[0,T^*).
\end{cases}
\ee
\end{lemma}
\pf
The proof of (\ref{non-energy inequality})$_1$ can be referred to \cite{Sun-Wang-Zhang 1} (Lemma 2). (\ref{non-energy inequality})$_2$ can be obtained by integrating (\ref{full N-S+1})$_1$ over $\mathbb{R}^3\times[0,t]$.
\endpf
\begin{lemma}\label{non-le:2.2}
Under the conditions of Theorem \ref{non-th:1.1} and (\ref{non-2.1}), if $3\mu>\lambda$, it holds that \be\label{u 2 nabla u 2} \sup\limits_{0\leq t\leq
T}\int_{\mathbb{R}^3}\rho|u|^4+\int_0^T\int_{\mathbb{R}^3}|u|^2|\nabla u|^2\,dx\leq C, \ee for any $T\in[0,T^*)$.
\end{lemma}
\pf  The proof of the lemma is quite similar to that of Lemma \ref{blow-up:le:5.1} except that $r=4$ and $P=\rho\theta$ here. From (\ref{blow-up:5.2}), we have
\beq\label{non-blow-up:5.2}\begin{split}
&\frac{d}{dt}\int_{\mathbb{R}^3}
\rho|u|^r+\int_{{\mathbb{R}^3}\cap\{|u|>0\}} r|u|^{r-2}\left(\mu|\nabla
u|^2+(\lambda+\mu)|\mathrm{div}u|^2+\mu(r-2)|\nabla|u||^2\right)\\
=&r\int_{{\mathbb{R}^3}\cap\{|u|>0\}}
\mathrm{div}(|u|^{r-2}u)P-r(r-2)(\mu+\lambda)\int_{{\mathbb{R}^3}\cap\{|u|>0\}}\mathrm{div}u|u|^{r-3}u\cdot\nabla|u|.
\end{split}\eeq
For any given $\varepsilon_1\in (0, 1)$, we define a function as in the proof of Lemma \ref{blow-up:le:5.1} as follows:
$$
\phi(\varepsilon_1,r)=\left\{\begin{array}{l}
\frac{\mu\varepsilon_1 (r-1)}{3\left(-\frac{4\mu}{3}-\lambda+
\frac{r^2(\mu+\lambda)}{4(r-1)}\right)}, \ \ \ {\rm if} \ \ \ \frac{r^2(\mu+\lambda)}{4(r-1)}-\frac{4\mu}{3}-\lambda>0, \\
[3mm] 0, \ \ \ {\rm otherwise}.
\end{array}
\right.
$$
{\bf Case 1:}
\be\label{non-blow-up:5.8}\int_{{\mathbb{R}^3}\cap\{|u|>0\}}|u|^r\left|\nabla\left(\frac{u}{|u|}\right)\right|^2>\phi(\varepsilon_1,r)\int_{{\mathbb{R}^3}\cap\{|u|>0\}}|u|^{r-2}\big|\nabla|u|\big|^2.\ee

Using the similar arguments like in the proof of Lemma \ref{blow-up:le:5.1}, for any $\varepsilon_0\in(0,1)$, we have \bex
\begin{split}&\frac{d}{dt}\int_{\mathbb{R}^3}
\rho|u|^r+\int_{{\mathbb{R}^3}\cap\{|u|>0\}} \mu r|u|^{r-2}\big|\nabla
|u|\big|^2+\int_{{\mathbb{R}^3}\cap\{|u|>0\}} \mu r|u|^r\left|\nabla
\left(\frac{u}{|u|}\right)\right|^2
\\&+\mu(r-2)r\int_{{\mathbb{R}^3}\cap\{|u|>0\}}|u|^{r-2}\big|\nabla|u|\big|^2
\\
\le& C\int_{{\mathbb{R}^3}\cap\{|u|>0\}}
\rho|u|^{r-2}\big|\nabla
|u|\big|+\mu r\varepsilon_0\int_{{\mathbb{R}^3}\cap\{|u|>0\}}
|u|^r\left|\nabla \left(\frac{u}{|u|}\right)\right|^2\\ &+\frac{C}{4\mu
r\varepsilon_0}\left(\int_{\mathbb{R}^3} \rho |u|^r\right)^\frac{r-2}{r}\left(\int_{\mathbb{R}^3}
\rho^{\frac{r}{2}+1}\right)^\frac{2}{r}
+\frac{r(r-2)^2(\mu+\lambda)}{4}\int_{{\mathbb{R}^3}\cap\{|u|>0\}}|u|^{r-2}\big|\nabla|u|\big|^2.
\end{split}\eex
Combining (\ref{non-2.1}) and (\ref{non-blow-up:5.8}), we have \beq\label{non-blow-up:5.10}
\begin{split}&\frac{d}{dt}\int_{\mathbb{R}^3} \rho|u|^r+rf(\varepsilon_0,\varepsilon_1, \varepsilon_2, r)\int_{{\mathbb{R}^3}\cap\{|u|>0\}}
|u|^{r-2}\big|\nabla |u|\big|^2+\mu r(1-\varepsilon_0)\varepsilon_2\int_{{\mathbb{R}^3}\cap\{|u|>0\}}|u|^r\left|\nabla
\left(\frac{u}{|u|}\right)\right|^2
\\ \le&C\int_{{\mathbb{R}^3}\cap\{|u|>0\}} \rho
|u|^{r-2}\big|\nabla |u|\big|+\frac{C}{4\mu
r\varepsilon_0}\left(\int_{\mathbb{R}^3} \rho
|u|^r\right)^\frac{r-2}{r},\end{split} \eeq
where \be\label{non-blow-up:5.12}
f(\varepsilon_0,\varepsilon_1, \varepsilon_2, r)=\mu
(1-\varepsilon_0)(1-\varepsilon_2)\phi(\varepsilon_1,r)+\mu
(r-1)-\frac{(r-2)^2(\mu+\lambda)}{4},
 \ee for $\varepsilon_2\in(0,1)$ decided later.

 ({\bf Sub-Case
1$_1$}):
 If $4\in\{r|\frac{r^2(\mu+\lambda)}{4(r-1)}-\frac{4\mu}{3}-\lambda>0\}$, i.e.,
 $\lambda>0$, we have
 \be\label{non-blow-up:5.11}\phi(\varepsilon_1,4)=
\frac{3\mu\varepsilon_1}{\lambda}.\ee
Substituting (\ref{non-blow-up:5.11}) into (\ref{non-blow-up:5.12}), we have \bex
f(\varepsilon_0,\varepsilon_1,\varepsilon_2,r)=\frac{\mu^2\varepsilon_1(1-\varepsilon_0)(1-\varepsilon_2)
(r-1)}{3\left(-\frac{4\mu}{3}-\lambda+
\frac{r^2(\mu+\lambda)}{4(r-1)}\right)}+\mu
(r-1)-\frac{(r-2)^2(\mu+\lambda)}{4}.
 \eex For $(\varepsilon_0,\varepsilon_1,\varepsilon_2,r)=(0,1,0,4)$, we have
\bex
f(0,1,0,4)=\frac{3\mu^2}{\lambda}+2\mu-\lambda>0,
\eex where we have used $0<\lambda<3\mu$.

Since $ f(\varepsilon_0,\varepsilon_1,\varepsilon_2,4)$ is continuous
 w.r.t. $(\varepsilon_0, \varepsilon_1,\varepsilon_2)$ over $[0, 1]\times [0, 1]\times
[0,1]$, there exist $\varepsilon_0, \varepsilon_1,\varepsilon_2 \in (0, 1)$
 such that
$$
f(\varepsilon_0,\varepsilon_1,\varepsilon_2,4)>0.
$$
By (\ref{non-blow-up:5.10}), Cauchy inequality and H\"older
inequality, for $r=4$, we have \bex\begin{split}&\frac{d}{dt}\int_{\mathbb{R}^3}
\rho|u|^4+4f(\varepsilon_0,\varepsilon_1,\varepsilon_2,4)\int_{{\mathbb{R}^3}\cap\{|u|>0\}}
|u|^2\big|\nabla |u|\big|^2+4\mu(1-\varepsilon_0)\varepsilon_2\int_{{\mathbb{R}^3}\cap\{|u|>0\}}|u|^4\left|\nabla
\left(\frac{u}{|u|}\right)\right|^2
\\ \le&2f(\varepsilon_0,\varepsilon_1,\varepsilon_2,4)\int_{{\mathbb{R}^3}\cap\{|u|>0\}} |u|^2\big|\nabla
|u|\big|^2+\frac{C}{f(\varepsilon_0,\varepsilon_1,\varepsilon_2,4)}\left(\int_{\mathbb{R}^3}
\rho |u|^4\right)^\frac{1}{2}\left(\int_{\mathbb{R}^3}
\rho^{3}\right)^\frac{1}{2}\\&+\frac{C}{\varepsilon_0}\left(\int_{\mathbb{R}^3} \rho
|u|^4\right)^\frac{1}{2}. \end{split}\eex This together with (\ref{non-2.1}) gives
\beq\label{non-blow-up:5.14}
\begin{split}
&\frac{d}{dt}\int_{\mathbb{R}^3} \rho|u|^4+2f(\varepsilon_0,\varepsilon_1,\varepsilon_2,4)\int_{{\mathbb{R}^3}\cap\{|u|>0\}}
|u|^2\big|\nabla |u|\big|^2+4\mu(1-\varepsilon_0)\varepsilon_2\int_{{\mathbb{R}^3}\cap\{|u|>0\}}|u|^4\left|\nabla
\left(\frac{u}{|u|}\right)\right|^2\\ \le&
C\left[\frac{1}{f(\varepsilon_0,\varepsilon_1,\varepsilon_2,4)}+\frac{1}{
\varepsilon_0}\right]\left(\int_{\mathbb{R}^3} \rho
|u|^4\right)^\frac{1}{2}.
\end{split}
\eeq
 ({\bf Sub-Case 1$_2$}): if $4\not\in\{r|\frac{r^2(\mu+\lambda)}{4(r-1)}-\frac{4\mu}{3}-\lambda>0\}$, i.e.,
 $\lambda\le0$, we have $\phi(\varepsilon_1,4)=0$.

In this case, it is easy to get \beq\label{non-blow-up:5.15}
\begin{split}&4f(\varepsilon_0,\varepsilon_1, \varepsilon_2, 4)=4(2\mu-\lambda)\ge 8\mu.\end{split}
 \eeq
By (\ref{non-blow-up:5.10}) (for $r=4$), (\ref{non-blow-up:5.15}), Cauchy inequality
and H\"older inequality, we have
 \bex
\begin{split}&\frac{d}{dt}\int_{\mathbb{R}^3}
\rho|u|^4+8\mu\int_{{\mathbb{R}^3}\cap\{|u|>0\}}
|u|^2\big|\nabla |u|\big|^2+4\mu (1-\varepsilon_0)\varepsilon_2\int_{{\mathbb{R}^3}\cap\{|u|>0\}}|u|^4\left|\nabla
\left(\frac{u}{|u|}\right)\right|^2 \\ \le&
C\int_{{\mathbb{R}^3}\cap\{|u|>0\}}\rho |u|^2\big|\nabla |u|\big|+\frac{C}{\varepsilon_0}\left(\int_{\mathbb{R}^3} \rho |u|^4\right)^\frac{1}{2}\\
\le&4\mu\int_{{\mathbb{R}^3}\cap\{|u|>0\}}
|u|^2\big|\nabla |u|\big|^2 +C\left(\int_{\mathbb{R}^3} \rho
|u|^4\right)^\frac{1}{2}\left(\int_{\mathbb{R}^3}
\rho^{3}\right)^\frac{1}{2}+\frac{C}{\varepsilon_0}\left(\int_{\mathbb{R}^3} \rho |u|^4\right)^\frac{1}{2}. \end{split}\eex Therefore,
\be\label{non-blow-up:5.16}\begin{split} &\frac{d}{dt}\int_{\mathbb{R}^3} \rho|u|^4+4\mu\int_{{\mathbb{R}^3}\cap\{|u|>0\}}
|u|^2\big|\nabla |u|\big|^2+4\mu (1-\varepsilon_0)\varepsilon_2\int_{{\mathbb{R}^3}\cap\{|u|>0\}}|u|^4\left|\nabla
\left(\frac{u}{|u|}\right)\right|^2\\ \le&
C\left(\int_{\mathbb{R}^3} \rho |u|^4\right)^\frac{1}{2}, \end{split}\ee where we have used (\ref{non-2.1}).

{\bf Case 2:} if
\be\label{non-blow-up:5.3}\int_{{\mathbb{R}^3}\cap\{|u|>0\}}|u|^r\left|\nabla\left(\frac{u}{|u|}\right)\right|^2\le\phi(\varepsilon_1,r)\int_{{\mathbb{R}^3}\cap\{|u|>0\}}|u|^{r-2}\big|\nabla|u|\big|^2.\ee

Using the similar arguments like in the proof of Lemma \ref{blow-up:le:5.1}, we have
 \bex
\begin{split}&\frac{d}{dt}\int_{\mathbb{R}^3} \rho|u|^r+\left[3r\left(\frac{4\mu}{3}+\lambda-
\frac{r^2(\mu+\lambda)}{4(r-1)}\right)\phi(\varepsilon_1,r)+\mu
r(r-1)\right]\int_{{\mathbb{R}^3}\cap\{|u|>0\}}
|u|^{r-2}\big|\nabla|u|\big|^2\\ \le&
C\int_{{\mathbb{R}^3}\cap\{|u|>0\}}
\rho^{1-\frac{r-2}{2r}}\rho^{\frac{r-2}{2r}}|u|^{r-2}|\nabla u|\\
\le&\varepsilon\int_{{\mathbb{R}^3}\cap\{|u|>0\}}|u|^{r-2}|\nabla
u|^2+\frac{C}{\varepsilon}\left(\int_{{\mathbb{R}^3}\cap\{|u|>0\}}
\rho|u|^r\right)^\frac{r-2}{r}\left(\int_{{\mathbb{R}^3}\cap\{|u|>0\}}
\rho^{\frac{r}{2}+1}\right)^\frac{2}{r}\\
\le&\varepsilon\big(1+\phi(\varepsilon_1,r)\big)\int_{{\mathbb{R}^3}\cap\{|u|>0\}}
|u|^{r-2}\big|\nabla|u|\big|^2+\frac{C}{\varepsilon}\left(\int_{{\mathbb{R}^3}\cap\{|u|>0\}}
\rho|u|^r\right)^\frac{r-2}{r},
\end{split}\eex
where we have used Cauchy inequality, H\"older inequality and (\ref{non-2.1}).

Taking
$\varepsilon=\big(2+2\phi(\varepsilon_1,r)\big)^{-1}\left[3r\left(\frac{4\mu}{3}+\lambda-
\frac{r^2(\mu+\lambda)}{4(r-1)}\right)\phi(\varepsilon_1,r)+\mu
r(r-1)\right]$, and using (\ref{non-blow-up:5.3}) and (\ref{blow-up:5.4}), we have \beq\label{non-blow-up:5.7}
\begin{split}&\frac{d}{dt}\int_{\mathbb{R}^3} \rho|u|^r+\frac{3r\left(\frac{4\mu}{3}+\lambda-
\frac{r^2(\mu+\lambda)}{4(r-1)}\right)\phi(\varepsilon_1,r)+\mu
r(r-1)}{2\left(1+\phi(\varepsilon_1,r)\right)}\int_{{\mathbb{R}^3}\cap\{|u|>0\}}
|u|^{r-2}\big|\nabla u\big|^2
\\ \le&\frac{C\big(1+\phi(\varepsilon_1,r)\big)}{3r\left(\frac{4\mu}{3}+\lambda-
\frac{r^2(\mu+\lambda)}{4(r-1)}\right)\phi(\varepsilon_1,r)+\mu
r(r-1)}\left(\int_{{\mathbb{R}^3}} \rho|u|^r\right)^\frac{r-2}{r},
\end{split}\eeq
for $r=4$.

By (\ref{blow-up:5.4}), (\ref{non-blow-up:5.14}), (\ref{non-blow-up:5.16}), (\ref{non-blow-up:5.7}) and Cauchy inequality, for {\bf Case
1} and {\bf Case 2}, we conclude that if $3\mu>\lambda$, there exists a constant $c_1>0$ such that \be\label{non-blow-up:5.17}
\frac{d}{dt}\int_{\mathbb{R}^3}\rho|u|^4+c_1\int_{{\mathbb{R}^3}\cap\{|u|>0\}}
|u|^2\big|\nabla u\big|^2\le C\int_{\mathbb{R}^3} \rho
|u|^4+C,\ee for $t\in[0,T^*)$. By
(\ref{non-blow-up:5.17}) and
Gronwall inequality, we get (\ref{u 2 nabla u 2}).
\endpf
\begin{lemma}\label{non-le: int nabla u}Under the conditions of Theorem \ref{non-th:1.1} and (\ref{non-2.1}), it holds that for any $T\in[0,T^*)$
\be\label{H 1 of u}\sup\limits_{0\le t\le T}\int_{\mathbb{R}^3}|\nabla u|^2\, dx+\int_0^T\int_{\mathbb{R}^3}\rho |u_t|^2\,dxdt\le C.\ee
\end{lemma}
\pf
Multiplying (\ref{full N-S+1})$_2$ by $u_t$, and integrating by parts over $\mathbb{R}^3$, we have
\beq\label{dt nabla u 2-1}\begin{split}
&\int_{\mathbb{R}^3}\rho |u_t|^2+\frac{1}{2}\frac{d}{dt}\int_{\mathbb{R}^3}\left(\mu|\nabla u|^2+(\mu+\lambda)|\mathrm{div}u|^2\right)\\=&-\int_{\mathbb{R}^3}\rho u\cdot\nabla u \cdot u_t+\frac{d}{dt}\int_{\mathbb{R}^3}P\mathrm{div}u-\int_{\mathbb{R}^3}P_t\mathrm{div}u\\=&\frac{d}{dt}\int_{\mathbb{R}^3}P\mathrm{div}u-\frac{1}{2(2\mu+\lambda)}
\frac{d}{dt}\int_{\mathbb{R}^3}P^2
-\frac{1}{2\mu+\lambda}\int_{\mathbb{R}^3}P_tG-\int_{\mathbb{R}^3}\rho u\cdot\nabla u \cdot u_t\\=&\sum\limits_{i=1}^4III_i,
\end{split}
\eeq where $G=(2\mu+\lambda)\mathrm{div}u-P$.

For $III_3$, recalling $\rho E=P+\frac{\rho|u|^2}{2}$, we have
\beq\label{I 3}\begin{split}
III_3=&-\frac{1}{2\mu+\lambda}\int_{\mathbb{R}^3}(\rho E)_t G+
\frac{1}{2\mu+\lambda}\int_{\mathbb{R}^3}\left(\frac{\rho|u|^2}{2}\right)_t G\\=&\sum\limits_{i=1}^2III_{3,i}.
\end{split}
\eeq
For $III_{3,1}$, using (\ref{full N-S})$_3$, integration by parts, (\ref{non-2.1}) and (\ref{non-energy inequality}), we have
\beq\label{I 3,1}\begin{split}
III_{3,1}=&-\frac{1}{2\mu+\lambda}\int_{\mathbb{R}^3}\rho \theta u\cdot\nabla G-\frac{1}{2\mu+\lambda}\int_{\mathbb{R}^3}\rho\frac{|u|^2}{2}u\cdot\nabla G-\frac{1}{2\mu+\lambda}\int_{\mathbb{R}^3}P
u\cdot\nabla G\\&+\frac{1}{2\mu+\lambda}\int_{\mathbb{R}^3}\mathcal {T}
u\nabla G+\frac{1}{2\mu+\lambda}\int_{\mathbb{R}^3}\nabla\theta\cdot\nabla G\\ \le&-\frac{1}{2\mu+\lambda}\int_{\mathbb{R}^3}\rho\frac{|u|^2}{2}u\cdot\nabla G+C\|\nabla G\|_{L^2}\left(\|\rho \theta u\|_{L^2}+\|Pu\|_{L^2}+\big\|u|\nabla u|\big\|_{L^2}+\|\nabla\theta\|_{L^2}\right)\\ \le&-\frac{1}{2\mu+\lambda}\int_{\mathbb{R}^3}\rho\frac{|u|^2}{2}u\cdot\nabla G+C\|\nabla G\|_{L^2}\left(\big\|u|\nabla u|\big\|_{L^2}+\|\nabla\theta\|_{L^2}+1\right).
\end{split}
\eeq
Taking $\mathrm{div}$ and $\mathrm{curl}$ on both side of (\ref{full N-S+1})$_2$, we get
\be\label{equation of G}
\Delta G=\mathrm{div}(\rho u_t+\rho u\cdot\nabla u),
\ee
and
\be\label{equation of curlu}
\mu\Delta(\mathrm{curl}u)=\mathrm{curl}(\rho u_t+\rho u\cdot\nabla u).
\ee
From the standard elliptic estimates together with (\ref{non-2.1}), we get
\beq\label{H 1 of G}
\begin{split}\|\nabla G\|_{L^{2}}\les \|\rho u_t\|_{L^{2}}+\|\rho u\cdot\nabla u\|_{L^2}\les
 \|\sqrt{\rho} u_t\|_{L^{2}}+\big\|\sqrt{\rho} |u||\nabla u|\big\|_{L^2},
\end{split}
\eeq
and
\beq\label{H 1 of curl u}
\begin{split}\|\nabla \mathrm{curl}u\|_{L^{2}} \les \|\rho u_t\|_{L^{2}}+\|\rho u\cdot\nabla u\|_{L^2}\les
 \|\sqrt{\rho} u_t\|_{L^{2}}+\big\|\sqrt{\rho} |u||\nabla u|\big\|_{L^2}.
\end{split}
\eeq
To handle the second term of the right hand side of (\ref{H 1 of G}) and (\ref{H 1 of curl u}), we use the fact
\be\label{equation of f}
-\Delta f=\nabla\times(\mathrm{curl}f)-\nabla\mathrm{div}f,\ \mathrm{in}\ \mathbb{R}^3,
\ee for some $f: \mathbb{R}^3\rightarrow\mathbb{R}^3$.
Using (\ref{equation of f}) and the elliptic estimates, we have
\be\label{nabla div curl}
\|\nabla f\|_{L^p}\les \|\mathrm{curl}f\|_{L^p}+\|\mathrm{div}f\|_{L^p},
\ee for any $p\in(1, \infty)$.
Let's go back to handle $\big\|\sqrt{\rho} |u||\nabla u|\big\|_{L^2}$. Using H\"older inequality, (\ref{non-2.1}), (\ref{u 2 nabla u 2}), (\ref{nabla div curl}) for $p=4$, Gagliardo-Nirenberg inequality, (\ref{H 1 of G}), (\ref{H 1 of curl u}) and Cauchy inequality, we have
\bex\begin{split}
\big\|\sqrt{\rho} |u||\nabla u|\big\|_{L^2}\les&\|\rho^\frac{1}{4} u\|_{L^4}\|\nabla u\|_{L^4}\les\|\mathrm{curl}u\|_{L^4}+\|\mathrm{div} u\|_{L^4}\\ \les&
\|\mathrm{curl}u\|_{L^4}+\|G\|_{L^4}+1\\ \le&C
\|\mathrm{curl}u\|_{L^2}^\frac{1}{4}\|\nabla\mathrm{curl}u\|_{L^2}^\frac{3}{4}+C\|G\|_{L^2}^\frac{1}{4}\|\nabla G\|_{L^2}^\frac{3}{4}+C\\ \le&
C\|\mathrm{curl}u\|_{L^2}^\frac{1}{4}\|\sqrt{\rho} u_t\|_{L^2}^\frac{3}{4}+C\|G\|_{L^2}^\frac{1}{4}\|\sqrt{\rho} u_t\|_{L^2}^\frac{3}{4}+\frac{1}{2}\big\|\sqrt{\rho} |u||\nabla u|\big\|_{L^2}+C\|\nabla u\|_{L^2}+C.
\end{split}
\eex
This, together with Young inequality, gives
\be\label{rho u nabla u}\begin{split}
\big\|\sqrt{\rho} |u||\nabla u|\big\|_{L^2} \le\epsilon\|\sqrt{\rho} u_t\|_{L^2}+C_\epsilon\|\nabla u\|_{L^2}+C,
\end{split}
\ee for any $\epsilon>0$.
Substituting (\ref{rho u nabla u}) into (\ref{H 1 of G}), we have
\beq\label{H 1 of G+1}
\begin{split}\|\nabla G\|_{L^{2}}\les
 \|\sqrt{\rho} u_t\|_{L^{2}}+\|\nabla u\|_{L^2}+1.
\end{split}
\eeq
Substituting (\ref{H 1 of G+1}) into (\ref{I 3,1}), and using Cauchy inequality, we have
\beq\label{I 3,1+1}\begin{split}
III_{3,1}\le&-\frac{1}{2\mu+\lambda}\int_{\mathbb{R}^3}\rho\frac{|u|^2}{2}u\cdot\nabla G+\frac{1}{6}\|\sqrt{\rho} u_t\|_{L^{2}}^2+
C\big\|u|\nabla u|\big\|_{L^2}^2\\&+C\|\nabla\theta\|_{L^2}^2+C\|\nabla u\|_{L^2}^2+C.
\end{split}
\eeq
For $III_{3,2}$, we have
\be\label{I 3,2}\begin{split}
III_{3,2}=&\frac{1}{2\mu+\lambda}\int_{\mathbb{R}^3}\frac{\rho_t|u|^2}{2} G+
\frac{1}{2\mu+\lambda}\int_{\mathbb{R}^3}\rho u\cdot u_t G\\ \le&-\frac{1}{2\mu+\lambda}\int_{\mathbb{R}^3}\frac{\mathrm{div}(\rho u)|u|^2}{2} G+
\frac{1}{24}\int_{\mathbb{R}^3}\rho |u_t|^2+C\int_{\mathbb{R}^3}\rho |u|^2 |G|^2\\ \le&\frac{1}{2\mu+\lambda}\int_{\mathbb{R}^3}\rho u\cdot\nabla u\cdot u G+\frac{1}{2\mu+\lambda}\int_{\mathbb{R}^3}\frac{\rho u|u|^2}{2}\cdot\nabla G+
\frac{1}{24}\int_{\mathbb{R}^3}\rho |u_t|^2\\&+C\int_{\mathbb{R}^3}\rho |u|^2 |\nabla u|^2+C\\ \le&C\int_{\mathbb{R}^3}\rho |u|^2|\nabla u|^2+\frac{1}{2\mu+\lambda}\int_{\mathbb{R}^3}\frac{\rho u|u|^2}{2}\cdot\nabla G+
\frac{1}{12}\int_{\mathbb{R}^3}\rho |u_t|^2+C.
\end{split}
\ee
Using (\ref{rho u nabla u}) again (for $\epsilon>0$ sufficiently small), together with (\ref{I 3,2}), Lemma \ref{non-le:2.2}, (\ref{H 1 of G+1}) and Cauchy inequality, we get
\be\label{I 3,2+1}\begin{split}
III_{3,2}\le&\frac{1}{12}\int_{\mathbb{R}^3}\rho |u_t|^2+\frac{1}{2\mu+\lambda}\int_{\mathbb{R}^3}\frac{\rho u|u|^2}{2}\cdot\nabla G+
\frac{1}{12}\int_{\mathbb{R}^3}\rho |u_t|^2+C\int_{\mathbb{R}^3} |\nabla u|^2+C\\ =&\frac{1}{6}\int_{\mathbb{R}^3}\rho |u_t|^2+\frac{1}{2\mu+\lambda}\int_{\mathbb{R}^3}\frac{\rho u|u|^2}{2}\cdot\nabla G+C\int_{\mathbb{R}^3} |\nabla u|^2+C.
\end{split}
\ee
Substituting (\ref{I 3,1+1}) and (\ref{I 3,2+1}) into (\ref{I 3}), we have
\beq\label{I 3+1}\begin{split}
III_3\le \frac{1}{3}\int_{\mathbb{R}^3}\rho |u_t|^2+
C\big\|u|\nabla u|\big\|_{L^2}^2+C\|\nabla\theta\|_{L^2}^2+C\|\nabla u\|_{L^2}^2+C.
\end{split}
\eeq
For $III_4$, using Cauchy inequality and (\ref{rho u nabla u}) (for $\epsilon>0$ sufficiently small), we have
\be\label{I 4}\begin{split}
III_4\le&\frac{1}{12}\int_{\mathbb{R}^3}\rho |u_t|^2+C\int_{\mathbb{R}^3}\rho |u|^2 |\nabla u|^2\\ \le&\frac{1}{6}\int_{\mathbb{R}^3}\rho |u_t|^2+C\int_{\mathbb{R}^3} |\nabla u|^2+C.
\end{split}
\ee
Putting (\ref{I 3+1}) and (\ref{I 4}) into (\ref{dt nabla u 2-1}), and integrating it over $[0,t]$, for $t<T^*$, we have
\bex\begin{split}
&\int_0^t\int_{\mathbb{R}^3}\rho |u_t|^2+\int_{\mathbb{R}^3}\left(\mu|\nabla u|^2+(\mu+\lambda)|\mathrm{div}u|^2\right)\\ \le&2\int_{\mathbb{R}^3}P\mathrm{div}u+C\int_0^t\big\|u|\nabla u|\big\|_{L^2}^2+C\int_0^t\|\nabla\theta\|_{L^2}^2+C\int_0^t\|\nabla u\|_{L^2}^2+C\\ \le&(\mu+\lambda)\int_{\mathbb{R}^3}|\mathrm{div}u|^2+C,
\end{split}
\eex where we have used Cauchy inequality, (\ref{non-2.1}), Lemmas \ref{non-le:2.1} and \ref{non-le:2.2}.

Therefore,
\bex\begin{split}
\int_0^t\int_{\mathbb{R}^3}\rho |u_t|^2+\int_{\mathbb{R}^3}|\nabla u|^2 \le C,
\end{split}
\eex for $t\in[0,T^*)$.

\endpf

\begin{lemma}\label{le: int rho u t}Under the conditions of Theorem \ref{non-th:1.1} and (\ref{non-2.1}), it holds that for any $T\in[0,T^*)$
\be\label{H 2 of u}\sup\limits_{0\le t\le T}\int_{\mathbb{R}^3}(|\nabla \theta|^2+\rho|\u|^2)+\int_0^T\int_{\mathbb{R}^3}(\rho |\dot{\theta}|^2+|\nabla\u|^2)\le C.\ee
\end{lemma}
\pf
Using the similar arguments as (\ref{i-dt rho u t-1}), we obtain
\beq\label{dt rho u t-1}
\begin{split}
&\frac{1}{2}\frac{d}{dt}\int_{\mathbb{R}^3}\rho|\u|^2+\int_{\mathbb{R}^3}\left(\mu|\nabla\u|^2+(\mu+\lambda)|\di\u|^2\right)\\
=&\int_{\mathbb{R}^3}\left( P_t\di\u + u\otimes\nabla P:\nabla\u\right)+\mu\int_{\mathbb{R}^3}\Big(\di(\de u\otimes u)-\de(u\cdot\nabla u)\Big)\cdot\u \\
+&(\mu+\lambda)\int_{\mathbb{R}^3}\Big(\di(\nabla\di u\otimes u)-\nabla\di(u\cdot\nabla u)\Big)\cdot\u
=\sum\ls_{i=1}^3IV_i.
\end{split}
\eeq
For $IV_1$, using (\ref{full N-S+1})$_3$ and integration by parts, we have
\beq\label{II 1}
\begin{split}
IV_1=&\int_{\mathbb{R}^3}\Big( (\rho\theta)_t\di\u-P(\nabla u)^t:\nabla\u -\rho\theta u\cdot\nabla\di\u \Big)\\
=&\int_{\mathbb{R}^3}\Big((\rho\theta)_t\di\u +\mathrm{div}(\rho\theta u)\di\u-P(\nabla u)^t:\nabla\u  \Big) \\
=&\int_{\mathbb{R}^3}\Big(\rho\dot{\theta}\di\u-P(\nabla u)^t:\nabla\u  \Big)
\\ \les&\|\sqrt{\rho}\|_{L^\infty}\|\sqrt{\rho}\dot{\theta}\|_{L^2}\|\di\u\|_{L^2}+\|P\|_{L^4}\|\nabla u\|_{L^4}\|\nabla \u\|_{L^2}.
\end{split}
\eeq
For $IV_2$ and $IV_3$, by (\ref{i-II 2}) and (\ref{i-II 3}), we have
\beq\label{II 2}
\begin{split}
IV_2
\les\|\nabla \u\|_{L^2}\|\nabla u\|_{L^4}^2,
\end{split}
\eeq
and
\beq\label{II 3}
\begin{split}
IV_3
\les\|\nabla \u\|_{L^2}\|\nabla u\|_{L^4}^2.
\end{split}
\eeq
Substituting (\ref{II 1}), (\ref{II 2}) and (\ref{II 3}) into (\ref{dt rho u t-1}), and using Cauchy inequality and (\ref{non-2.1}), we have
\bex
\begin{split}
\frac{1}{2}\frac{d}{dt}\int_{\mathbb{R}^3}\rho|\u|^2+\int_{\mathbb{R}^3}\left(\mu|\nabla\u|^2+(\mu+\lambda)|\di\u|^2\right)
\le\frac{\mu}{2}\|\nabla \u\|_{L^2}^2+C\|\sqrt{\rho}\dot{\theta}\|_{L^2}^2+C\|\nabla u\|_{L^4}^4+C.
\end{split}
\eex
Integrating this inequality over $[0,t]$ for $t\in(0,T^*)$, we have
\beq\label{dt rho u t-2}
\begin{split}
\int_{\mathbb{R}^3}\rho|\u|^2+\int_0^t\int_{\mathbb{R}^3}|\nabla\u|^2\le C\int_0^t\|\sqrt{\rho}\dot{\theta}\|_{L^2}^2+C\int_0^t\|\nabla u\|_{L^4}^4+C.
\end{split}
\eeq
The next step is to get some estimates for $\theta$. We rewrite (\ref{full N-S+1})$_3$ as follows:
\be\label{rho theta+ =}
\rho \dot{\theta}+\rho\theta\mathrm{div}u=\frac{\mu}{2}\left|\nabla u+(\nabla u)^\prime\right|^2+\lambda(\mathrm{div}u)^2+\Delta\theta.
\ee
Multiplying (\ref{rho theta+ =}) by $\dot{\theta}$, and integrating by parts over $\mathbb{R}^3$, we have
\be\label{dt nabla theta}\begin{split}
\int_{\mathbb{R}^3}\rho |\dot{\theta}|^2+\frac{1}{2}\frac{d}{dt}\int_{\mathbb{R}^3}|\nabla \theta|^2=&-\int_{\mathbb{R}^3}\rho\theta\mathrm{div}u\dot{\theta}+\int_{\mathbb{R}^3}\left(\frac{\mu}{2}\left|\nabla u+(\nabla u)^\prime\right|^2+\lambda(\mathrm{div}u)^2\right)\theta_t\\&+\int_{\mathbb{R}^3}\left(\frac{\mu}{2}\left|\nabla u+(\nabla u)^\prime\right|^2+\lambda(\mathrm{div}u)^2\right)u\cdot\nabla\theta+\int_{\mathbb{R}^3}\Delta\theta u\cdot\nabla\theta\\=& \sum\limits_{i=1}^4V_i.
\end{split}
\ee
For $V_1$, using Cauchy inequality, (\ref{non-2.1}) and (\ref{H 1 of u}), we have
\be\label{III 1}\begin{split}
V_1\le\frac{1}{8}\int_{\mathbb{R}^3}\rho|\dot{\theta}|^2+C.
\end{split}
\ee For $V_2$, using H\"older inequality, (\ref{non-2.1}) and
(\ref{H 1 of u}), we have \bex\begin{split}
V_2=&\frac{d}{dt}\int_{\mathbb{R}^3}\left(\frac{\mu}{2}\left|\nabla
u+(\nabla
u)^\prime\right|^2+\lambda(\mathrm{div}u)^2\right)\theta-\mu\int_{\mathbb{R}^3}\left(\nabla
u+(\nabla u)^\prime\right):\left(\nabla u_t+(\nabla
u_t)^\prime\right)\theta\\&-2\lambda\int_{\mathbb{R}^3}\mathrm{div}u\mathrm{div}u_t\theta\\=&
\frac{d}{dt}\int_{\mathbb{R}^3}\left(\frac{\mu}{2}\left|\nabla
u+(\nabla u)^\prime\right|^2+\lambda(\mathrm{div}u)^2\right)\theta
-\mu\int_{\mathbb{R}^3}\left(\nabla u+(\nabla
u)^\prime\right):\left(\nabla \u+(\nabla
\u)^\prime\right)\theta\\&+\mu\int_{\mathbb{R}^3}\left(\nabla
u+(\nabla u)^\prime\right):\left(\nabla (u\cdot\nabla u)+(\nabla
(u\cdot\nabla
u))^\prime\right)\theta-2\lambda\int_{\mathbb{R}^3}\mathrm{div}u\mathrm{div}\u\theta\\&
+2\lambda\int_{\mathbb{R}^3}\mathrm{div}u\mathrm{div}(u\cdot\nabla
u)\theta\\ \le&
\frac{d}{dt}\int_{\mathbb{R}^3}\left(\frac{\mu}{2}\left|\nabla
u+(\nabla u)^\prime\right|^2+\lambda(\mathrm{div}u)^2\right)\theta
+C\|\nabla u\|_{L^2}\|\nabla
\u\|_{L^2}\\&+\mu\int_{\mathbb{R}^3}\left(\nabla u+(\nabla
u)^\prime\right):\left(\nabla (u\cdot\nabla u)+(\nabla
(u\cdot\nabla u))^\prime\right)\theta
+2\lambda\int_{\mathbb{R}^3}\mathrm{div}u\mathrm{div}(u\cdot\nabla
u)\theta.
\end{split}
\eex Using integration by parts, (\ref{H 1 of u}) and
(\ref{non-2.1}), we have \be\label{III 2-1}\begin{split} V_2\le&
\frac{d}{dt}\int_{\mathbb{R}^3}\left(\frac{\mu}{2}\left|\nabla
u+(\nabla u)^\prime\right|^2+\lambda(\mathrm{div}u)^2\right)\theta
+C\|\nabla \u\|_{L^2}\\&+\mu\int_{\mathbb{R}^3}\left(\nabla
u+(\nabla u)^\prime\right):\left(\nabla u\cdot\nabla u+(\nabla
u\cdot\nabla
u)^\prime\right)\theta\\&+\mu\int_{\mathbb{R}^3}\left(\nabla
u+(\nabla u)^\prime\right):u\cdot\nabla\left(\nabla u+(\nabla
u)^\prime\right)\theta
+2\lambda\int_{\mathbb{R}^3}\mathrm{div}u(\nabla u)^\prime:\nabla
u\theta\\&+2\lambda\int_{\mathbb{R}^3}u\cdot\nabla
\mathrm{div}u\mathrm{div}u\theta
\\ \le&
\frac{d}{dt}\int_{\mathbb{R}^3}\left(\frac{\mu}{2}\left|\nabla
u+(\nabla u)^\prime\right|^2+\lambda(\mathrm{div}u)^2\right)\theta
+C\|\nabla \u\|_{L^2}+C\int_{\mathbb{R}^3}|\nabla
u|^3\\&-\mu\int_{\mathbb{R}^3}\frac{|\nabla u+(\nabla
u)^\prime|^2}{2}\mathrm{div}u\theta-\mu\int_{\mathbb{R}^3}\frac{|\nabla
u+(\nabla
u)^\prime|^2}{2}u\cdot\nabla\theta\\&-\lambda\int_{\mathbb{R}^3}
(\mathrm{div}u)^3\theta-\lambda\int_{\mathbb{R}^3}|\mathrm{div}u|^2u\cdot\nabla\theta
\\ \le&
\frac{d}{dt}\int_{\mathbb{R}^3}\left(\frac{\mu}{2}\left|\nabla u+(\nabla u)^\prime\right|^2+\lambda(\mathrm{div}u)^2\right)\theta
+C\|\nabla \u\|_{L^2}+C\int_{\mathbb{R}^3}|\nabla u|^3+C\int_{\mathbb{R}^3}|\nabla u|^2|u||\nabla\theta|.
\end{split}
\ee
Using H\"older inequality, Cauchy inequality, Gagliardo-Nirenberg inequality and (\ref{H 1 of u}), we have
\be\label{non-nabla u u nabla theta}\begin{split}
\int_{\mathbb{R}^3}|\nabla u|^2|u||\nabla\theta|\les&\|\nabla u\|_{L^4}^2\|u\|_{L^6}\|\nabla\theta\|_{L^3}\\ \les&
\|\nabla u\|_{L^4}^4+\|\nabla u\|_{L^2}^2\|\nabla\theta\|_{L^2}\|\nabla^2\theta\|_{L^2}\le C\|\nabla u\|_{L^4}^4+C\|\nabla\theta\|_{L^2}\|\nabla^2\theta\|_{L^2}.
\end{split}
\ee
From the standard elliptic estimates and (\ref{rho theta+ =}), we have
\be\label{H 2 of theta}
\|\nabla^2\theta\|_{L^2}\les\|\rho \dot{\theta}\|_{L^2}+\|\rho\theta\mathrm{div}u\|_{L^2}+\|\nabla u\|_{L^4}^2\le C\|\sqrt{\rho}\dot{\theta}\|_{L^2}+C\|\nabla u\|_{L^4}^2+C,
\ee where we have used (\ref{non-2.1}) and (\ref{H 1 of u}).

Substituting (\ref{H 2 of theta}) into (\ref{non-nabla u u nabla theta}), and using Cauchy inequality, we have
\be\label{non-nabla u u nabla theta+1}\begin{split}
\int_{\mathbb{R}^3}|\nabla u|^2|u||\nabla\theta|\le \frac{1}{8}\|\sqrt{\rho}\dot{\theta}\|_{L^2}^2+C\|\nabla u\|_{L^4}^4+C\|\nabla\theta\|_{L^2}^2+C.
\end{split}
\ee Substituting (\ref{non-nabla u u nabla theta+1}) into (\ref{III 2-1}), we have
\be\label{III 2}\begin{split}
V_2\le&
\frac{d}{dt}\int_{\mathbb{R}^3}\left(\frac{\mu}{2}\left|\nabla u+(\nabla u)^\prime\right|^2+\lambda(\mathrm{div}u)^2\right)\theta
+C\|\nabla \u\|_{L^2}+C\int_{\mathbb{R}^3}|\nabla u|^3\\&+\frac{1}{8}\|\sqrt{\rho}\dot{\theta}\|_{L^2}^2+C\|\nabla u\|_{L^4}^4+C\|\nabla\theta\|_{L^2}^2+C.
\end{split}
\ee
For $V_3$, using (\ref{non-nabla u u nabla theta+1}), we have
\be\label{III 3}\begin{split}
V_3\les&\int_{\mathbb{R}^3}|\nabla u|^2|u||\nabla\theta|\le\frac{1}{8}\|\sqrt{\rho}\dot{\theta}\|_{L^2}^2+C\|\nabla u\|_{L^4}^4+C\|\nabla\theta\|_{L^2}^2+C.
\end{split}
\ee
For $V_4$, using H\"older inequality, Gagliardo-Nirenberg inequality, (\ref{H 1 of u}), (\ref{H 2 of theta}) and Young inequality, we have
\be\label{III 4}\begin{split}
V_4\les& \|\Delta\theta\|_{L^2}\|u\|_{L^6}\|\nabla\theta\|_{L^3}\les \|\Delta\theta\|_{L^2}\|\nabla u\|_{L^2}\|\nabla\theta\|_{L^2}^\frac{1}{2}\|\nabla^2\theta\|_{L^2}^\frac{1}{2}\\ \les&
\|\nabla\theta\|_{L^2}^\frac{1}{2}\|\nabla^2\theta\|_{L^2}^\frac{3}{2}\le \frac{1}{8}\|\sqrt{\rho}\dot{\theta}\|_{L^2}^2+C\|\nabla\theta\|_{L^2}^2
+C\|\nabla u\|_{L^4}^4+C.
\end{split}
\ee
Putting (\ref{III 1}), (\ref{III 2}), (\ref{III 3}) and (\ref{III 4}) into (\ref{dt nabla theta}), and integrating the resulting inequality over $[0,t]$ for $t\in(0,T^*)$, we have
\be\label{dt nabla theta+1}\begin{split}
\int_0^t\int_{\mathbb{R}^3}\rho |\dot{\theta}|^2+\int_{\mathbb{R}^3}|\nabla \theta|^2\le&
2\int_{\mathbb{R}^3}\left(\frac{\mu}{2}\left|\nabla u+(\nabla u)^\prime\right|^2+\lambda(\mathrm{div}u)^2\right)\theta+C\int_0^t\|\nabla \u\|_{L^2}\\&+C
\int_0^t\int_{\mathbb{R}^3}|\nabla u|^3+
C\int_0^t\|\nabla u\|_{L^4}^4+C\\ \le&C\int_0^t\|\nabla \u\|_{L^2}+C\int_0^t\int_{\mathbb{R}^3}|\nabla u|^3+
C\int_0^t\|\nabla u\|_{L^4}^4+C,
\end{split}
\ee where we have used (\ref{non-2.1}), (\ref{non-energy inequality}) and (\ref{H 1 of u}). Multiplying (\ref{dt nabla theta+1}) by $2C$, and adding the resulting inequality into (\ref{dt rho u t-2}), we have
\bex
\begin{split}
&C\int_0^t\int_{\mathbb{R}^3}\rho |\dot{\theta}|^2+2C\int_{\mathbb{R}^3}|\nabla \theta|^2+\int_{\mathbb{R}^3}\rho|\u|^2+\int_0^t\int_{\mathbb{R}^3}|\nabla\u|^2\\ \le&2C^2\int_0^t\|\nabla \u\|_{L^2}+2C^2\int_0^t\int_{\mathbb{R}^3}|\nabla u|^3+2C^2\int_0^t\|\nabla u\|_{L^4}^4+C\int_0^t\|\nabla u\|_{L^4}^4+C.
\end{split}
\eex
This together with Cauchy inequality, we have
\be\label{non-sum}
\begin{split}
&\int_{\mathbb{R}^3}(|\nabla \theta|^2+\rho|\u|^2)+\int_0^t\int_{\mathbb{R}^3}(\rho |\dot{\theta}|^2+|\nabla\u|^2) \les\int_0^t\int_{\mathbb{R}^3}|\nabla u|^3+\int_0^t\|\nabla u\|_{L^4}^4+1\\ &\les\int_0^t\int_{\mathbb{R}^3}|\mathrm{curl} u|^3+\int_0^t\int_{\mathbb{R}^3}|G|^3+\int_0^t\|\mathrm{curl}u\|_{L^4}^4+\int_0^t\|G\|_{L^4}^4
+1\\ &\les\int_0^t\|\mathrm{curl} u\|_{L^2}^\frac{3}{2}\|\nabla\mathrm{curl}u\|_{L^2}^\frac{3}{2}+\int_0^t\|G\|_{L^2}^\frac{3}{2}\|\nabla G\|_{L^2}^\frac{3}{2}+\int_0^t\|\mathrm{curl}u\|_{L^2}\|\nabla\mathrm{curl}u\|_{L^2}^3\\&+\int_0^t\|G\|_{L^2}\|\nabla G\|_{L^2}^3
+1,
\end{split}
\ee where we have used (\ref{nabla div curl}) and Gagliardo-Nirenberg inequality. By (\ref{non-sum}), (\ref{non-2.1}), (\ref{non-energy inequality}), (\ref{H 1 of u}), (\ref{H 1 of curl u}), (\ref{rho u nabla u}) and (\ref{H 1 of G+1}), we have
\be\label{non-sum+1}
\begin{split}
\int_{\mathbb{R}^3}(|\nabla \theta|^2+\rho|\u|^2)+\int_0^t\int_{\mathbb{R}^3}(\rho |\dot{\theta}|^2+|\nabla\u|^2)\les\int_0^t\|\sqrt{\rho} u_t\|_{L^2}^3+1.
\end{split}
\ee
From (\ref{rho u nabla u}), we have
\bex
\begin{split}
\|\sqrt{\rho}u_t\|_{L^2}\le\|\sqrt{\rho}\u\|_{L^2}+\|\sqrt{\rho}u\cdot\nabla u\|_{L^2}\le\|\sqrt{\rho}\u\|_{L^2}+\epsilon\|\sqrt{\rho} u_t\|_{L^2}+C_\epsilon\|\nabla u\|_{L^2}+C.
\end{split}
\eex
Taking $\epsilon=\frac{1}{2}$, using (\ref{H 1 of u}), we have
\be\label{rho u t and rho u.}
\begin{split}
\|\sqrt{\rho}u_t\|_{L^2}\les\|\sqrt{\rho}\u\|_{L^2}+1.
\end{split}
\ee
Substituting (\ref{rho u t and rho u.}) into (\ref{non-sum+1}), and using Cauchy inequality and (\ref{H 1 of u}), we have
\be\label{non-sum+2}
\begin{split}
\int_{\mathbb{R}^3}(|\nabla \theta|^2+\rho|\u|^2)+\int_0^t\int_{\mathbb{R}^3}(\rho |\dot{\theta}|^2+|\nabla\u|^2)\les\int_0^t(\|\sqrt{\rho} u_t\|_{L^2}\|\sqrt{\rho} \u\|_{L^2}^2)+1.
\end{split}
\ee
Since $\|\sqrt{\rho} u_t\|_{L^2}$ is bounded in $L^1-$norm over $(0, t)$ (see (\ref{H 1 of u})), we use (\ref{non-sum+2}) and Gronwall inequality to get (\ref{H 2 of u}).
\endpf
\begin{corollary}Under the conditions of Theorem \ref{non-th:1.1} and (\ref{non-2.1}), it holds that for any $T\in[0,T^*)$
\be\label{u infty}\sup\limits_{0\le t\le T}\left(\|\nabla G\|_{L^2}+\|\nabla \mathrm{curl}u\|_{L^2}+\|\nabla u\|_{L^6}+\|u\|_{L^\infty}\right)+\int_0^T(\|\mathrm{div}u\|_{L^\infty}^2+\|\nabla^2\theta\|_{L^2}^2)\le C.\ee
\end{corollary}
\pf
It follows from (\ref{equation of G}) and (\ref{equation of curlu}), we have
\be\label{L 2 of nabla G}
\|\nabla G\|_{L^2}\les \|\rho\u\|_{L^2}\le C,
\ee
and
\be\label{L infty of divu}\begin{split}
\int_0^T\|\mathrm{div}u\|_{L^\infty}^2\les& \int_0^T\|G\|_{L^\infty}^2+1\les \int_0^T\|G\|_{L^6}^2+\int_0^T\|\nabla G\|_{L^6}^2+1\\ \les& \int_0^T\|\nabla G\|_{L^2}^2+\int_0^T\|\rho\u\|_{L^6}^2+1\les \int_0^T\|\nabla\u\|_{L^2}^2+1\le C,
\end{split}\ee
and
\be\label{L 2 of nabla curl u}
\|\nabla \mathrm{curl}u\|_{L^2}\les \|\rho\u\|_{L^2}\le C,
\ee
where we have used (\ref{non-2.1}), (\ref{H 2 of u}) and Sobolev inequality.

By (\ref{nabla div curl}), we have
\be\label{L 6 of nabla u}\begin{split}
\|\nabla u\|_{L^6}\les& \|\mathrm{div} u\|_{L^6}+\|\mathrm{curl} u\|_{L^6}\les\|G\|_{L^6}+\|\mathrm{curl} u\|_{L^6}+1\\ \les&\|\nabla G\|_{L^2}+\|\nabla\mathrm{curl} u\|_{L^2}+1\le C,
\end{split}
\ee where we have used (\ref{non-2.1}), (\ref{non-energy inequality}), Sobolev inequality, (\ref{L 2 of nabla G}) and (\ref{L 2 of nabla curl u}).

By (\ref{H 1 of u}), (\ref{L 6 of nabla u}) and Sobolev inequality, we have
$$
\|u\|_{L^\infty}\les \|u\|_{L^6}+\|\nabla u\|_{L^6}\les \|\nabla u\|_{L^2}+\|\nabla u\|_{L^6}\le C.
$$
Using (\ref{H 2 of theta}), (\ref{H 2 of u}), the interpolation inequality, (\ref{H 1 of u}) and (\ref{L 6 of nabla u}), we get
\bex\begin{split}
\int_0^T\int_{\mathbb{R}^3}|\nabla^2\theta|^2\les& \int_0^T\int_{\mathbb{R}^3}\rho|\dot{\theta}|^2+\int_0^T\|\nabla u\|_{L^4}^4+1\\ \les&
\int_0^T\|\nabla u\|_{L^2}\|\nabla u\|_{L^6}^3+1\le C.
\end{split}
\eex
\endpf

\begin{lemma}\label{le: H 2 of theta}Under the conditions of Theorem \ref{non-th:1.1} and (\ref{non-2.1}), it holds that for any $T\in[0,T^*)$
\be\label{non-H 2 of theta}
\sup\limits_{0\le t\le T}\int_{\mathbb{R}^3}\rho|\theta_t|^2+\int_0^T\int_{\mathbb{R}^3}|\nabla\theta_t|^2\le C.
\ee
\end{lemma}
\pf Differentiating (\ref{full N-S+1})$_3$ with respect to $t$, multiplying it by $\theta_t$, and using integration by parts, we have
\be\label{dt H 2 of theta}\begin{split}
&\frac{1}{2}\frac{d}{dt}\int_{\mathbb{R}^3}\rho|\theta_t|^2+\int_{\mathbb{R}^3}|\nabla\theta_t|^2\\ \le&-
\int_{\mathbb{R}^3}\rho_t\left(\frac{\theta_t}{2}+u\cdot\nabla\theta+\theta\mathrm{div}u\right)\theta_t-\int_{\mathbb{R}^3}\rho( u_t\cdot\nabla\theta+u\cdot\nabla\theta_t+\theta_t\mathrm{div}u)\theta_t-\int_{\mathbb{R}^3}\rho\theta\mathrm{div}u_t\theta_t\\&+\mu\int_{\mathbb{R}^3}\left(\nabla u+(\nabla u)^\prime\right):\left(\nabla u_t+(\nabla u_t)^\prime\right)\theta_t+2\lambda\int_{\mathbb{R}^3}\mathrm{div}u\mathrm{div}u_t\theta_t=\sum\limits_{i=1}^5VI_i.
\end{split}
\ee
For $VI_1$, we have
\be\label{IV 1}\begin{split}
VI_1=&\int_{\mathbb{R}^3}\mathrm{div}(\rho u)\left(\frac{\theta_t}{2}+u\cdot\nabla\theta+\theta\mathrm{div}u\right)\theta_t\\=&
-\int_{\mathbb{R}^3}\rho u\cdot\nabla\theta_t\left(\frac{\theta_t}{2}+u\cdot\nabla\theta+\theta\mathrm{div}u\right)-\int_{\mathbb{R}^3}\rho u\cdot\frac{\nabla\theta_t}{2}\theta_t\\&-\int_{\mathbb{R}^3}\rho u\cdot\left(\nabla u\cdot\nabla\theta+u\cdot\nabla\nabla\theta\right)\theta_t-\int_{\mathbb{R}^3}\rho u\cdot\left(\nabla\theta\mathrm{div}u+\theta\nabla\mathrm{div}u\right)\theta_t\\=&\sum\limits_{i=1}^4VI_{1,i}.
\end{split}
\ee
For $VI_{1,1}$, we have
\be\label{IV 1,1}\begin{split}
VI_{1,1} \le&
\frac{1}{24}\int_{\mathbb{R}^3}|\nabla\theta_t|^2+C\int_{\mathbb{R}^3}\rho^2|u|^2|\theta_t|^2+C\int_{\mathbb{R}^3}\rho^2|u|^4|\nabla\theta|^2
+C\int_{\mathbb{R}^3}\rho^2|u|^2|\theta|^2|\mathrm{div}u|^2\\ \le&
\frac{1}{24}\int_{\mathbb{R}^3}|\nabla\theta_t|^2+C\int_{\mathbb{R}^3}\rho|\theta_t|^2
+C,
\end{split}
\ee where we have used Cauchy inequality, (\ref{non-2.1}), (\ref{H 1 of u}), (\ref{H 2 of u}) and (\ref{u infty}).

For $VI_{1,2}$, using Cauchy inequality, (\ref{non-2.1}) and (\ref{u infty}) again, we have
\be\label{IV 1,2}\begin{split}
VI_{1,2}\le&\frac{1}{24}\int_{\mathbb{R}^3}|\nabla\theta_t|^2+C\int_{\mathbb{R}^3}\rho |\theta_t|^2.
\end{split}
\ee
For $VI_{1,3}$, using Cauchy inequality, (\ref{non-2.1}) and (\ref{u infty}) again, along with H\"older inequality, Gagliardo-Nirenberg inequality and (\ref{H 2 of u}), we have
\be\label{IV 1,3}\begin{split}
VI_{1,3}\les&\int_{\mathbb{R}^3}\rho |\theta_t|^2+
\int_{\mathbb{R}^3}|\nabla u|^2|\nabla\theta|^2+\int_{\mathbb{R}^3}|\nabla\nabla\theta|^2\\ \les&\int_{\mathbb{R}^3}\rho |\theta_t|^2+
\|\nabla u\|_{L^6}^2\|\nabla\theta\|_{L^3}^2+\int_{\mathbb{R}^3}|\nabla\nabla\theta|^2\\ \les&\int_{\mathbb{R}^3}\rho |\theta_t|^2+\|\nabla\theta\|_{L^2}\|\nabla^2\theta\|_{L^2}+\int_{\mathbb{R}^3}|\nabla\nabla\theta|^2\\ \les&\int_{\mathbb{R}^3}\rho |\theta_t|^2+\int_{\mathbb{R}^3}|\nabla^2\theta|^2+1.
\end{split}
\ee
For $VI_{1,4}$, we have
\be\label{IV 1,4}\begin{split}
VI_{1,4}=&-\int_{\mathbb{R}^3}\rho u\cdot\nabla\theta\mathrm{div}u\theta_t-\int_{\mathbb{R}^3}\rho\theta u\cdot\nabla\mathrm{div}u\theta_t\\ \les&
\int_{\mathbb{R}^3}\rho |\theta_t|^2+\int_{\mathbb{R}^3}|\nabla\theta|^2|\mathrm{div}u|^2-\frac{1}{2\mu+\lambda}\int_{\mathbb{R}^3}\rho\theta u\cdot\nabla G\theta_t-\frac{1}{2\mu+\lambda}\int_{\mathbb{R}^3}\rho\theta u\cdot\nabla (\rho\theta)\theta_t\\ \les&
\int_{\mathbb{R}^3}\rho |\theta_t|^2+\|\nabla\theta\|_{L^3}^2\|\mathrm{div}u\|_{L^6}^2+\int_{\mathbb{R}^3}|\nabla G|^2-\frac{1}{2\mu+\lambda}\int_{\mathbb{R}^3}\rho^2\theta u\cdot\nabla\theta\theta_t\\&-\frac{1}{2\mu+\lambda}\int_{\mathbb{R}^3}\rho\theta^2 u\cdot\nabla \rho\theta_t\\ \le&
C\int_{\mathbb{R}^3}\rho |\theta_t|^2+C\|\nabla^2\theta\|_{L^2}+C+\frac{1}{2\mu+\lambda}\int_{\mathbb{R}^3}\frac{\rho^2}{2}\theta^2 \mathrm{div}u \theta_t+\frac{1}{2\mu+\lambda}\int_{\mathbb{R}^3}\frac{\rho^2}{2}\theta^2 u\cdot\nabla\theta_t
\\&+\frac{1}{2\mu+\lambda}\int_{\mathbb{R}^3}\rho^2\theta u\cdot\nabla\theta\theta_t \le
C\int_{\mathbb{R}^3}\rho |\theta_t|^2+C\|\nabla^2\theta\|_{L^2}^2+\frac{1}{24}\int_{\mathbb{R}^3}|\nabla\theta_t|^2+C,
\end{split}
\ee where we have used Cauchy inequality, (\ref{non-2.1}), (\ref{u infty}), H\"older inequality, Gagliardo-Nirenberg inequality, (\ref{H 2 of u}), integration by parts, (\ref{non-energy inequality}) and (\ref{H 1 of u}).

Substituting (\ref{IV 1,1}), (\ref{IV 1,2}), (\ref{IV 1,3}) and (\ref{IV 1,4}) into (\ref{IV 1}), we have
\be\label{IV 1+1}\begin{split}
VI_1\le&\frac{1}{8}\int_{\mathbb{R}^3}|\nabla\theta_t|^2
+C\int_{\mathbb{R}^3}\rho |\theta_t|^2+C\int_{\mathbb{R}^3}|\nabla^2\theta|^2+C.
\end{split}
\ee
For $VI_2$, using Cauchy inequality, H\"older inequality, (\ref{non-2.1}) and (\ref{u infty}), we have
\bex\begin{split}
VI_2=&-\int_{\mathbb{R}^3}\rho u_t\cdot\nabla\theta\theta_t-\int_{\mathbb{R}^3}\rho u\cdot\nabla\theta_t\theta_t-\int_{\mathbb{R}^3}\rho|\theta_t|^2\mathrm{div}u\\ \le&-\int_{\mathbb{R}^3}\rho \u\cdot\nabla\theta\theta_t+\int_{\mathbb{R}^3}\rho (u\cdot\nabla)u\cdot\nabla\theta\theta_t+\frac{1}{8}\int_{\mathbb{R}^3}|\nabla\theta_t|^2+C\left(\|\mathrm{div}u\|_{L^\infty}+1\right)\int_{\mathbb{R}^3}\rho|\theta_t|^2
\\ \le&C\|\sqrt{\rho}\theta_t\|_{L^2}\|\u\|_{L^6}\|\nabla\theta\|_{L^3}+C\|\sqrt{\rho}\theta_t\|_{L^2}\|\nabla u\|_{L^6}\|\nabla\theta\|_{L^3}+\frac{1}{8}\int_{\mathbb{R}^3}|\nabla\theta_t|^2\\&+C\left(\|\mathrm{div}u\|_{L^\infty}+1\right)\int_{\mathbb{R}^3}\rho|\theta_t|^2.
\end{split}
\eex This together with Sobolev inequality, Gagliardo-Nirenberg inequality, (\ref{u infty}) and (\ref{H 2 of u}), we have
\be\label{IV 2}\begin{split}
VI_2 \le&C\|\sqrt{\rho}\theta_t\|_{L^2}\left(\|\nabla\u\|_{L^2}+1\right)\|\nabla\theta\|_{L^2}^\frac{1}{2}\|\nabla^2\theta\|_{L^2}^\frac{1}{2}
+\frac{1}{8}\int_{\mathbb{R}^3}|\nabla\theta_t|^2+C\left(\|\mathrm{div}u\|_{L^\infty}+1\right)\int_{\mathbb{R}^3}\rho|\theta_t|^2\\ \le&\frac{1}{8}\int_{\mathbb{R}^3}|\nabla\theta_t|^2+C\left(\|\mathrm{div}u\|_{L^\infty}+\|\nabla\u\|_{L^2}^2+1\right)\int_{\mathbb{R}^3}\rho|\theta_t|^2 +C\|\nabla^2\theta\|_{L^2}^2+C.
\end{split}
\ee
For $VI_3$, we have
\be\label{IV 3}\begin{split}
VI_3=&-\int_{\mathbb{R}^3}\rho\theta\mathrm{div}\u\theta_t+\int_{\mathbb{R}^3}\rho\theta\mathrm{div}(u\cdot\nabla u)\theta_t\\ \le&
C\int_{\mathbb{R}^3}\rho|\theta_t|^2+C\int_{\mathbb{R}^3}|\mathrm{div}\u|^2+C\int_{\mathbb{R}^3}\rho|\theta||\nabla u|^2|\theta_t|+\int_{\mathbb{R}^3}\rho\theta\theta_t u\cdot\nabla \mathrm{div}u\\ \le&C\int_{\mathbb{R}^3}\rho|\theta_t|^2+C\int_{\mathbb{R}^3}|\mathrm{div}\u|^2+C\int_{\mathbb{R}^3}|\nabla u|^4+\frac{1}{2\mu+\lambda}\int_{\mathbb{R}^3}\rho\theta\theta_t u\cdot\nabla G\\&+\frac{1}{2\mu+\lambda}\int_{\mathbb{R}^3}\rho^2\theta\theta_t u\cdot\nabla\theta+\frac{1}{2\mu+\lambda}\int_{\mathbb{R}^3}\rho\theta^2\theta_t u\cdot\nabla\rho\\ \le&C\int_{\mathbb{R}^3}\rho|\theta_t|^2+C\int_{\mathbb{R}^3}|\mathrm{div}\u|^2+\frac{1}{2\mu+\lambda}\int_{\mathbb{R}^3}\rho\theta^2\theta_t u\cdot\nabla\rho+C,
\end{split}
\ee where we have used (\ref{non-2.1}), Cauchy inequality, the interpolation inequality, (\ref{H 1 of u}), (\ref{H 2 of u}) and (\ref{u infty}).

 To handle the third term of the right hand side of (\ref{IV 3}), we use integration by parts. More precisely,
\be\label{the third term}
\begin{split}
\frac{1}{2\mu+\lambda}\int_{\mathbb{R}^3}\rho\theta^2\theta_t u\cdot\nabla\rho=&
-\frac{1}{2\mu+\lambda}\int_{\mathbb{R}^3}\frac{\rho^2}{2}\theta^2\theta_t \mathrm{div}u-\frac{1}{2\mu+\lambda}\int_{\mathbb{R}^3}\frac{\rho^2}{2}\theta^2 u\cdot\nabla\theta_t\\&-\frac{1}{2\mu+\lambda}\int_{\mathbb{R}^3}\rho^2\theta\theta_t u\cdot\nabla\theta\\ \le&C\int_{\mathbb{R}^3}\rho|\theta_t|^2+\frac{1}{8}\int_{\mathbb{R}^3}|\nabla\theta_t|^2+C,
\end{split}
\ee where we have used Cauchy inequality, (\ref{non-2.1}), (\ref{non-energy inequality}), (\ref{H 1 of u}), (\ref{H 2 of u}) and (\ref{u infty}).

Substituting (\ref{the third term}) into (\ref{IV 3}), we have
\be\label{IV 3+1}\begin{split}
VI_3 \le&\frac{1}{8}\int_{\mathbb{R}^3}|\nabla\theta_t|^2+C\int_{\mathbb{R}^3}\rho|\theta_t|^2+C\int_{\mathbb{R}^3}|\mathrm{div}\u|^2+C.
\end{split}
\ee
Similar to $V_2$, for $VI_4$ and $VI_5$, we deduce
\be\label{IV 4 and IV 5}
\begin{split}
VI_4+VI_5\le&C\|\nabla\u\|_{L^2}\|\nabla u\|_{L^3}\|\theta_t\|_{L^6}+C\int_{\mathbb{R}^3}|\nabla u|^3|\theta_t|+C\int_{\mathbb{R}^3}|\nabla u|^4+\frac{1}{16}\int_{\mathbb{R}^3}|\nabla\theta_t|^2\\ \le&\frac{1}{16}\int_{\mathbb{R}^3}|\nabla\theta_t|^2 +C\|\nabla\u\|_{L^2}\|\theta_t\|_{L^6}+C\|\nabla u\|_{L^\frac{18}{5}}^3\|\theta_t\|_{L^6}+C\\ \le&\frac{1}{16}\int_{\mathbb{R}^3}|\nabla\theta_t|^2 +C\left(\|\nabla\u\|_{L^2}+1\right)\|\nabla\theta_t\|_{L^2}+C
\\ \le&\frac{1}{8}\int_{\mathbb{R}^3}|\nabla\theta_t|^2+C\int_{\mathbb{R}^3}|\nabla\u|^2+C,
\end{split}
\ee
where we have used H\"older inequality, integration by parts, Cauchy inequality, (\ref{H 1 of u}), (\ref{u infty}), the interpolation inequality and Sobolev inequality.

Putting (\ref{IV 1+1}), (\ref{IV 2}), (\ref{IV 3+1}) and (\ref{IV 4 and IV 5}) into (\ref{dt H 2 of theta}), we have
\be\label{dt H 2 of theta+1}\begin{split}
\frac{d}{dt}\int_{\mathbb{R}^3}\rho|\theta_t|^2+\int_{\mathbb{R}^3}|\nabla\theta_t|^2\le&
C\left(\|\mathrm{div}u\|_{L^\infty}+\|\nabla\u\|_{L^2}^2+1\right)\int_{\mathbb{R}^3}\rho|\theta_t|^2 \\&+C\int_{\mathbb{R}^3}(|\nabla\u|^2+|\nabla^2\theta|^2)+C.
\end{split}
\ee
By (\ref{dt H 2 of theta+1}), (\ref{H 2 of u}), (\ref{u infty}) and Gronwall inequality, we complete the proof of Lemma \ref{le: H 2 of theta}.
\endpf
\begin{corollary}Under the conditions of Theorem \ref{non-th:1.1} and (\ref{non-2.1}), it holds that for any $T\in[0,T^*)$
\be\label{cor:H 2 of theta}
\sup\limits_{0\le t\le T}\int_{\mathbb{R}^3}|\nabla^2\theta|^2\le C.
\ee
\end{corollary}
\pf It follows from (\ref{H 2 of theta}), (\ref{non-2.1}), (\ref{H 1 of u}), (\ref{H 2 of u}), (\ref{u infty}), (\ref{non-H 2 of theta}) and the interpolation inequality that
\bex
\|\nabla^2\theta\|_{L^2}\le C\|\sqrt{\rho}\theta_t\|_{L^2}+C\|\sqrt{\rho}u\cdot\nabla\theta\|_{L^2}+C\le C.
\eex
\endpf
\begin{lemma}Under the conditions of Theorem \ref{non-th:1.1} and (\ref{non-2.1}), it holds that for any $T\in[0,T^*)$
\be\label{nabla rho}\sup\limits_{0\le t\le T}\left(\|\nabla\rho\|_{L^l}+\|\rho_t\|_{L^l}\right)\le C,\ee for $l=2,q$.
\end{lemma}
\pf The proof of the lemma is similar to the arguments as in \cite{Huang-Li-Xin: Serrin, Sun-Wang-Zhang}. We omit it for brevity.
\endpf
\begin{corollary}Under the conditions of Theorem \ref{non-th:1.1} and (\ref{non-2.1}), it holds that for any $T\in[0,T^*)$
\be\label{cor:H 2 of u}
\sup\limits_{0\le t\le T}\int_{\mathbb{R}^3}\left(\rho|u_t|^2+|\nabla^2u|^2\right)+\int_0^T\left(\|u_t\|_{D^1}^2+\|(u, \theta)\|_{D^{2,q}}^2\right)\le C.
\ee
\end{corollary}
\pf
Replacing $f$ in (\ref{equation of f}) by $u$, and using the elliptic estimates, (\ref{non-2.1}), (\ref{H 2 of u}), (\ref{u infty}) and (\ref{nabla rho}), we get
\be\label{H 2 of u+1}
\begin{split}
\|\nabla^2u\|_{L^2}\les& \|\nabla \mathrm{curl}u\|_{L^2}+\|\nabla\mathrm{div}u\|_{L^2}\les\|\nabla G\|_{L^2}+\|\nabla P(\rho,\theta)\|_{L^2}+1\\ \les&
\|\nabla\rho\|_{L^2}+\|\nabla\theta\|_{L^2}+1\le C.
\end{split}
\ee
It follows from (\ref{non-2.1}), (\ref{H 1 of u}), (\ref{H 2 of u}), (\ref{u infty}) and (\ref{H 2 of u+1}) that
\bex\begin{split}
\int_{\mathbb{R}^3}\rho|u_t|^2\les \int_{\mathbb{R}^3}\rho|\u|^2+\int_{\mathbb{R}^3}\rho|u\cdot\nabla u|^2\le C,
\end{split}\eex
and
\bex\begin{split}
\int_0^t\int_{\mathbb{R}^3}|\nabla u_t|^2\les \int_0^t\int_{\mathbb{R}^3}|\nabla \u|^2+\int_0^t\int_{\mathbb{R}^3}|\nabla (u\cdot\nabla u)|^2\le C.
\end{split}\eex
By (\ref{full N-S+1})$_2$, H\"older inequality, (\ref{non-2.1}), (\ref{H 2 of u}), Sobolev inequality, (\ref{cor:H 2 of theta}) and (\ref{nabla rho}), we get
\be\label{nabla 2 u q}\begin{split}
\int_0^t\|\nabla^2u\|_{L^q}^2\les&\int_0^t\|\rho\u\|_{L^q}^2+\int_0^t\|\nabla P(\rho,\theta)\|_{L^q}^2\les\int_0^t\|\u\|_{L^6}^2+\int_0^t\|\nabla\rho\|_{L^q}^2+\int_0^t\|\nabla\theta\|_{L^q}^2\\ \les&\int_0^t\|\nabla\u\|_{L^2}^2+\int_0^t\|\nabla\theta\|_{L^2}^2+\int_0^t\|\nabla^2\theta\|_{L^2}^2+1\le C.
\end{split}\ee
Using H\"older inequality, (\ref{non-2.1}), (\ref{H 2 of u}), Sobolev inequality and (\ref{cor:H 2 of theta}) again, together with (\ref{full N-S+1})$_3$, (\ref{H 1 of u}), (\ref{u infty}), (\ref{non-H 2 of theta}), (\ref{H 2 of u+1}) and (\ref{nabla 2 u q}), we get
\bex\begin{split}
\int_0^t\|\nabla^2\theta\|_{L^q}^2\les&\int_0^t\|\rho\theta_t\|_{L^q}^2+\int_0^t\|\rho u\cdot\nabla\theta\|_{L^q}^2+\int_0^t\|\rho\theta\mathrm{div}u\|_{L^q}^2+\int_0^t\left\||\nabla u|^2\right\|_{L^q}^2\\ \le&C\int_0^t\left\||\nabla u|^2\right\|_{L^q}^2+C\le C\int_0^t\|\nabla u\|_{L^\infty}^2\left\|\nabla u\right\|_{L^q}^2+C\\ \le&C\int_0^t\left\|\nabla^2 u\right\|_{L^q}^2+C\le C.
\end{split}\eex
\endpf \\

By (\ref{H 1 of u}), (\ref{H 2 of u}), (\ref{non-H 2 of theta}), (\ref{cor:H 2 of u}), (\ref{cor:H 2 of theta}) and (\ref{nabla rho}), we get (\ref{non-uniform_est1}) and (\ref{non-uniform_est2}). Thus, the proof of Theorem \ref{non-th:1.1} is complete. \endpf

\section*{Acknowledgements} This work was supported by the
National Natural Science Foundation of China $\#$10625105,
$\#$11071093, the PhD specialized grant of the Ministry of Education
of China $\#$20100144110001, and the Special Fund for Basic
Scientific Research  of Central Colleges $\#$CCNU10C01001.


\vskip 1cm


\addcontentsline{toc}{section}{\\References}

\begin{thebibliography}{99}
\bibitem{BKM}
J. T. Beale, T. Kato, A. Majda.
{\em Remarks on the breakdown of smooth solutions for the 3-D Euler equation}. Comm. Math. Phys., 94(1984), 61-66.

\bibitem{Bresch-Desjardins} D. Bresch, B. Desjardins. {\em On the existence of global weak solutions to the Navier-Stokes
equations for viscous compressible and heat conducting fluids.} J.
Math. Pures Appl., 87(2007), 57-90.

\bibitem{Carlson} J. Carlson, A. Jaffe, A. Wiles. {\em The Millennuim
Prize Problems.} The American Mathematical Society, Providence, RI,
2006.

\bibitem{Cho-Choe-Kim} Y. Cho, H. J. Choe, H. Kim. {\em Unique solvability of the initial boundary value problems for compressible viscous fluids}. J. Math. Pures Appl. 83(2004), 243-275.

\bibitem{Cho-Kim; perfect gas} Y. Cho, H. Kim. {\em Existence results for viscous polytropic fluids with vacuum.} J. Differential Equations, 228(2006), 377-411.

\bibitem{Cho-Kim: classical} Y. Cho, H. Kim. {\em On classical solutions of the compressible Navier-Stokes equations with
nonnegative initial densities.} Manuscripta Math., 120(2006),
91-129.

\bibitem{Choe-Bum} H.J. Choe, B.J. Jin. {\em Regularity of weak solutions of the compressible Navier-Stokes equations}. J. Korean
Math. Soc., 40 (2003), 1031-1050.

\bibitem{Choe-Kim:symmetric} H.J. Choe, H. Kim. {\em Global existence of the radially symmetric solutions of
the Navier-Stokes equations for the isentropic compressible fluids.}
Math. Methods Appl. Sci., 28(2005), 1-28.

\bibitem{Ding-Wen-Zhu} S.J. Ding, H.Y. Wen, C.J. Zhu. {\em Global classical large solutions to 1D compressible
Navier-Stokes equations with density-dependent viscosity and vacuum}. J. Differential Equations, 251(2011), 1696-1725.

\bibitem{Fan-Jiang} J.S. Fan, S. Jiang. {\em Blow-Up criteria for the navier-stokes equations of compressible fluids}. J.Hyper. Diff.
Eqs., 5(2008), 167-185.

\bibitem{Fan-Jiang-Ni} J.S. Fan, S. Jiang, G. Ni. {\em Uniform boundedness
of the radially symmetric solutions of the Navier-Stokes equations
for isentropic compressible fluids.} Osaka J. Math., 46(2009),
863-876.

\bibitem{Fan-Jiang-Ou} J.S. Fan, S. Jiang, Y.B. Ou. {\em A blow-up criterion for compressible
viscous heat-conductive flows}. Ann. I. H. Poincar¨¦-AN, 27(2010),
337-350.

\bibitem{Fang-Zi-Zhang} D.Y. Fang, R.Z. Zi, T. Zhang. {\em A blow-up criterion for two dimensional
compressible viscous heat-conductive flows}. arXiv:1107.4663v1 [math.AP] 23 Jul 2011.

\bibitem{Feireisl-book} E. Feireisl, {\em Dynamics of Viscous Compressible Fluids,} Oxford Univ. Press, Oxford, 2004.

\bibitem{Feireisl2} E. Feireisl, A. Novotn$\acute{\mathrm{y}}$, H.
Petzeltov$\mathrm{\acute{a}}$. {\em On the existence of globally
defined weak solutions to the Navier-Stokes equations.} J. Math.
Fluid Mech., 3(2001), 358-392.

\bibitem{Hoff1995}
D. Hoff. {\em Global solutions of the Navier-Stokes equations for
multidimensional compressible flow with discontinuous initial data}.
J. Differential Equations, 120(1995), 215-254.

\bibitem{Hoff ARMA} D. Hoff. {\em Discontinuous solutions of the Navier-Stokes
equations for multidimensional flows
of heat-conducting fluids}. Arch. Rational Mech. Anal., 139(1997), 303-354.

\bibitem{Huang-Wang-Wen} T. Huang, C.Y. Wang, H.Y. Wen. {\em Strong solutions of the compressible nematic liquid crystal flow}. J. Differential Equations (2011), doi:10.1016/j.jde.2011.07.036.

\bibitem{Huang-Wang-Wen:arma} T. Huang, C.Y. Wang, H.Y. Wen. Blow up criterion for compressible nematic liquid
crystal ows in dimension three, Arch. Rational Mech. Anal., 2011, to appear.

\bibitem{Huang-Li-Xin: Serrin}
X.D. Huang, J. Li, Z.P. Xin. {\em Serrin type criterion for the
three-dimensional viscous compressible flows}. SIAM J. Math. Anal., 43(2011), 1872-1886.

\bibitem{Huang-Li-Xin: classical} X.D. Huang, J. Li, Z.P. Xin. {\em Global well-posedness of classical solutions with large oscillations and vacuum to the three-dimensional isentropic compressible Navier-Stokes equations}. arXiv:1004.4749.

\bibitem{Huang-Li-Xin: blow up} X.D. Huang, J. Li, Z.P. Xin. {\em Blowup criterion for the compressible flows with vacuum
states}. Commun. Math. Phys., 301(2011), 23-35.



\bibitem{Itaya} N. Itaya. {\em On the Cauchy problem for the system of fundamental equations describing the movement of compressible viscous
fluid.} Kodai Math. Sem. Rep., 23(1971), 60-120.

\bibitem{Jiang} S. Jiang, P. Zhang. {\em On spherically symmetric solutions
of the compressible isentropic Navier-Stokes equations.} Comm. Math.
Phys., 215(2001), 559-581.

\bibitem{Jiang1} S. Jiang. {\em Global spherically symmetric solutions
to the equations of a viscous polytropic ideal gas in an exterior
domain.} Comm. Math. Phys., 178(1996), 339-374.

\bibitem{Kawohl} B. Kawohl. {\em Global existence of large solutions to initial boundary value
problems for a viscous, heat-conducting, one-dimensional real gas.}
J. Differential Equations, 58(1985), 76-103.

\bibitem{Kazhikhov-Shelukhi} A.V. Kazhikhov, V.V. Shelukhi. {\em Unique global solution with respect to time of
initial-boundary value problems for one-dimensional equations of a
viscous gas.} Prikl. Mat. Meh., 41(1977), 282-291.

\bibitem{Lions} P. L. Lions, Mathematical topics in fluid mechanics. Vol. 2. Compressible models. Oxford Lecture Series in Mathematics and its Applications, 10. Oxford Science Publications. The Clarendon Press, Oxford University Press, New York, 1998. xiv+348 pp.

\bibitem{Luo} T. Luo, Z.P. Xin, T. Yang. {\em Interface behavior of
compressible Navier-Stokes equations with vacuum.} SIAM J. Math.
Anal., 31(2000), 1175-1191.

\bibitem{Matsumura-Nishida: Kyoto Un} A. Matsumura, T. Nishida. {\em The initial value problem for the equations of motion of
viscous and heat-conductive gases.} J. Math. Kyoto Univ., 20(1980),
67-104.

\bibitem{Matsumura-Nishida: CMP} A. Matsumura, T. Nishida. {\em The initial boundary value problems for the equations of
motion of compressible and heat-conductive fluids.} Comm. Math.
Phys., 89(1983), 445-464.

\bibitem{Ponce} G. Ponce. {\em Remarks on a paper: "Remarks on the breakdown of smooth solutions for
the 3-D Euler equations"}.  Comm. Math. Phys., 98(1985), 349-353.

\bibitem{Rona} O. Rozanova. {\em Blow-up of smooth highly decreasing at infinity
solutions to the compressible Navier¨CStokes equation}. J. Differential Equations, 245(2008), 1762-1774.


\bibitem{salvi}
R. Salvi, I. Stra$\breve{\mathrm{s}}$kraba. {\em Global existence
for viscous compressible fluids and their behavior as
$t\rightarrow\infty$.} J. Fac. Sci. Univ. Tokyo, Sect. IA Math.,
40(1993), 17-51.

\bibitem{serrin} J. Serrin. {\em On the interior regularity of weak solutions of the Navier-Stokes equations}.
Arch. Rational Mech. Anal., 9(1962), 187-195.

\bibitem{Sun-Wang-Zhang}
Y.Z. Sun, C. Wang, Z.F. Zhang. {\em A Beale-Kato-Majda blow-up criterion
for the 3-D compressible Navier-Stokes equations}. J. Math. Pures
Appl., 95(2011), 36-47.

\bibitem{Sun-Wang-Zhang 1} Y.Z. Sun, C. Wang, Z.F. Zhang. {\em A Beale-Kato-Majda criterion for three
dimensional compressible viscous
heat-conductive Flows}. Arch. Rational Mech. Anal., Digital Object Identifier (DOI) 10.1007/s00205-011-0407-1.

\bibitem{Tani} A. Tani. {\em On the first initial-boundary value problem of compressible viscous fluid
motion.} Publ. Res. Inst. Math. Sci. Kyoto Univ., 13(1977), 193-253.

\bibitem{Wen-Yao-Zhu} H.Y. Wen, L. Yao, C.J. Zhu. {\em A blow-up criterion of strong solution to a 3D viscous liquid-gas
two-phase flow model with vacuum.} J. Math. Pures Appl. (2011), doi:10.1016/j.matpur.2011.09.005.

\bibitem{Wen-Zhu 1} H.Y. Wen, C.J. Zhu. {\em Global classical large solutions to Navier-Stokes
equations for viscous compressible and heat conducting fluids with
vacuum}. arXiv:1103.1421v1 [math.AP] 8 Mar 2011, 1-38.

\bibitem{Wen-Zhu 2} H.Y. Wen, C.J. Zhu. {\em Global symmetric classical and strong solutions of
the full compressible Navier-Stokes equations with
vacuum and large initial data}. arXiv:1109.5328v1 [math.AP] 25 Sep 2011.

\bibitem{Xin} Z.P. Xin. {\em Blowup of smooth solutions to the compressible Navier-Stokes equation with
compact density.} Comm. Pure Appl. Math., 51(1998), 229-240.

\end{thebibliography}
\end{document}